\documentclass[11pt, letterpaper]{article}
\pdfoutput=1 
\usepackage[authoryear]{natbib}
\usepackage{amsmath,amssymb,amsfonts,amsthm,bbm}
\usepackage{epic,eepic,epsfig,longtable}
\usepackage{multirow,verbatim}
\usepackage{array}

\usepackage{epsfig}
\usepackage{setspace}
\usepackage{CJK}
\usepackage{color}

\makeatletter
\def\singlespace{\def\baselinestretch{1}\@normalsize}

\makeatletter
\def\singlespace{\def\baselinestretch{1}\@normalsize}

\numberwithin{equation}{section}

\renewcommand{\hat}{\widehat}

\newcommand{\bfm}[1]{\ensuremath{\mathbf{#1}}}

          \def\cL{{\cal  L}}

          \def\cR{{\cal  R}}

\def\bx{\bfm x}

\newcommand{\bfsym}[1]{\ensuremath{\boldsymbol{#1}}}

\def\bbeta     {\bfsym \beta}

\def\bnu       {\bfsym {\nu}}

\def\bDelta  {\bfsym {\Delta}}

\renewcommand{\hat}{\widehat}

\DeclareMathOperator{\argmin}{argmin}

\DeclareMathOperator{\E}{E}

\def\newpage{\vfill\eject}

\def\var{\mbox{var}}

\def\today{\ifcase\month\or
  January\or February\or March\or April\or May\or June\or
  July\or August\or September\or October\or November\or December\fi
  \space\number\day, \number\year}

\newdimen\biblioindent    \biblioindent=30pt

 at 8truept

\def\la{\lambda}

\newcommand{\beq}{\begin{equation}}
  \newcommand{\eeq}{\end{equation}}
\newcommand{\beqn}{\begin{eqnarray}}
  \newcommand{\eeqn}{\end{eqnarray}}
\newcommand{\beqnn}{\begin{eqnarray*}}
  \newcommand{\eeqnn}{\end{eqnarray*}}

\def\bbC {\mathbb{C}}

\allowdisplaybreaks
\usepackage{verbatim}
\renewcommand{\baselinestretch}{1.66}
\newtheorem{lem}{Lemma}
\newtheorem{dfn}{Definition}
\newtheorem{thm}{Theorem}
\newcounter{CondCounter}

\newcommand{\lonenorm}[1]{\lVert#1\rVert_1}
\newcommand{\ltwonorm}[1]{\lVert#1\rVert_2}

\newcommand{\supnorm}[1]{\left \lVert#1 \right \rVert_{\infty}}

\def \betaa    {\bbeta_{\alpha}^{\ast}}
\def \betas    {\bbeta^{\ast}}

\def \betacs   {\bbeta_{\alpha}^{c\ast}}
\def \sqLogpN  {\sqrt{(\log p)/n}}
\def \la       {\ell_{\alpha}}
\def \kl       {\kappa_l}
\def \ku       {\kappa_u}
\def \ko       {\kappa_0}
\def \Rp       {\mathbb{R}^p}
\def \fd       {f_{\bDelta}}
\def \gd       {g_{\bDelta}}
\def \Rq       {R_{q}}
\def \bbP      {\mathbb{P}}
\def \bbS      {\mathbb{S}}
\def \sidx     {\sup_{\bDelta\in \bbS_2(1)\cap \bbS_1(t)}}

\begin{document}

\title{Robust Estimation of High-Dimensional Mean Regression\footnote{
Supported in part by NSF Grants DMS-1206464 and DMS-1406266 and NIH grants R01-GM072611-9 and NIH R01GM100474-4.}
}
\author{
  Jianqing Fan, Quefeng Li, Yuyan Wang \vspace{0.1in}\\ 
  Department of Operations Research and Financial Engineering,\\
  Princeton University, Princeton, NJ 08544}
\date{October 6, 2014}

\maketitle
\begin{abstract}
  Data subject to heavy-tailed errors are commonly encountered in various scientific fields,
  especially in the modern era with explosion of massive data.  To address this problem, procedures
  based on quantile regression and Least Absolute Deviation (LAD) regression have been developed in
  recent years. These methods essentially estimate the conditional median (or quantile) function.
  They can be very different from the conditional mean functions when distributions are asymmetric and heteroscedastic.  How can we efficiently estimate
  the mean regression functions in ultra-high dimensional setting with existence of only the second
  moment?  To solve this problem, we propose a penalized Huber loss with diverging parameter to
  reduce biases created by the traditional Huber loss.  Such a penalized robust approximate
  quadratic (RA-quadratic) loss will be called
  RA-Lasso.  In the ultra-high dimensional setting, where the dimensionality can grow exponentially
  with the sample size, our results reveal that the RA-lasso estimator produces a consistent
  estimator at the same rate as the optimal rate under the light-tail situation.  We further study
  the computational convergence of RA-Lasso and show that the composite gradient descent algorithm
  indeed produces a solution that admits the same optimal rate after sufficient iterations.  As a
  byproduct, we also establish the concentration inequality for estimating population mean when
  there exists only the second moment.  We compare RA-Lasso with other regularized robust
  estimators based on quantile regression and LAD regression. Extensive simulation studies
  demonstrate the satisfactory finite-sample performance of RA-Lasso.
\end{abstract}


\newpage

\section{Introduction}
\label{sec:introduction}
Our era has witnessed the massive explosion of data and a dramatic improvement of technology in
collecting and processing large data sets. We often encounter huge data sets that the number of
features greatly surpasses the number of observations. It makes many traditional statistical
analysis tools infeasible and poses great challenge on developing new tools. Regularization methods
have been widely used for the analysis of high-dimensional data. These methods penalize the least
squares or the likelihood function with the $L_1$-penalty on the unknown parameters (Lasso,
\cite{TIB96}), or a folded concave penalty function such as the SCAD \citep{FL01} and the
MCP\citep{ZHA10}. However, these penalized least-squares methods are sensitive to the tails of the
error distributions, particularly for ultrahigh dimensional covariates, as the maximum spurious
correlation between the covariates and the realized noises can be large in those cases.  As a
result, theoretical properties are often obtained under light-tailed error distributions
\citep{BRT09, FL11}.

To tackle the problem of heavy-tailed errors, robust regularization methods have been extensively
studied. \cite{LZ08}, \cite{WL09} and \cite{ZY08} developed robust regularized estimators
based on quantile regression for the case of fixed dimensionality. \cite{BC11} studied the
$L_1$-penalized quantile regression in high dimensional sparse models. \cite{FFB14} further
considered an adaptively weighted $L_1$ penalty to alleviate the bias problem and showed the oracle
property and asymptotic normality of the corresponding estimator. Other robust estimators were
developed based on Least Absolute Deviation (LAD) regression. \cite{WAN13} studied the
$L_1$-penalized LAD regression and showed that the estimator achieves near oracle risk performance
under the high dimensional setting.

The above methods essentially estimate the conditional
\emph{median (or quantile)} regression, instead of the conditional \emph{mean} regression function.
In the applications where the mean regression is of interest, these methods are not feasible
unless a strong assumption is made that the distribution of errors is symmetric around
zero. A simple example is the heteroscedastic linear model with asymmetric noise distribution.  Another example is to estimate the conditional variance function such as ARCH model \citep{ENG82}.  In these cases, the conditional mean and conditional median are very different.  Another important example is to estimate large covariance matrix without assuming light-tails.  We will explain this more in details in Section 5.  In addition, LAD-based methods tend to penalize strongly on small errors. If only a small proportion of samples are outliers, they are expected to be less efficient than the least squares based method.

A natural question is then how to conduct ultrahigh dimensional mean regression when the tails of
errors are not light?  How to estimate the sample mean with very fast concentration when the
distribution has only bounded second moment?  These simple questions have not been carefully
studied.  LAD-based methods do not intend to answer these questions as they alter the problems of
the study.  This leads us to consider Huber loss as another way of robustification. The Huber loss
\citep{HUB64} is a hybrid of squared loss for relatively small errors and absolute loss for
relatively large errors, where the degree of hybridization is controlled by one tuning
parameter. Unlike the traditional Huber loss, we allow the regularization parameter to diverge (or
converge if its reciprocal is used) in order to reduce the bias induced by the Huber loss for estimating conditional mean
regression function.  In this paper, we consider the regularized approximate quadratic (RA-Lasso)
estimator with an $L_1$ penalty and show that it admits the same $L_2$ error rate as the optimal
error rate in the light-tail situation. In particular, if the distribution of errors is indeed
symmetric around 0 (where the median and mean agree), this rate is the same as the regularized LAD
estimator obtained in \cite{WAN13}. Therefore, the RA-Lasso estimator does not lose efficiency in
this special case. In practice, since the distribution of errors is unknown, RA-Lasso is more
flexible than the existing methods in terms of estimating the conditional mean regression function.

A by-product of our method is that the RA-Lasso estimator of the population mean has the
exponential type of concentration even in presence of the finite second moment. \cite{CAT12}
studied this type of problem and proposed a class of losses to result in a robust $M$-estimator of
mean with exponential type of concentration. We further extend his idea to the sparse linear
regression setting.

As done in many other papers, estimators with nice sampling properties are typically defined
through the optimization of a target function such as the penalized least-squares.  The properties
that are established are not the same as the ones that are computed. Following the framework of
\cite{ANW12}, we propose the composite gradient descent algorithm for solving the RA-Lasso
estimator and develop the sampling properties by taking computational error into consideration. We
show that the algorithm indeed produces a solution that admits the same optimal $L_2$ error rate as
the theoretical estimator after sufficient number of iterations.

This paper is organized as follows. First, in Section \ref{sec2}, we introduce the RA-Lasso
estimator and show that it has the same $L_2$ error rate as the optimal rate under light-tails. In
Section \ref{sec3}, we study the property of the composite gradient descent algorithm for solving
our problem and show that the algorithm produces a solution that performs as well as the solution
in theory. In Section \ref{sec4}, we show the connection between Huber loss and Catoni loss and
establish an concentration inequality for robust estimation of mean.  The estimation of the error's
variance is investigated in Section \ref{sec5}. Numerical studies are given in Section \ref{sec6}
and \ref{sec7} to compare our method with two competitors. All technical proofs are presented in
Section \ref{sec:proofs}.

\section{RA-Lasso estimator}
\label{sec2}
We consider the linear regression model
\begin{equation}
  \label{eq2.1}
  y_i=\bx_i^T\bbeta^{\ast}+\epsilon_i,
\end{equation}
where $\{\bx_i\}_{i=1}^n$ are independent and identically distributed (i.i.d) $p$-dimensional
covariate vectors, $\{\epsilon_i\}_{i=1}^n$ are i.i.d errors, and $\bbeta^{\ast}$ is a
$p$-dimensional regression coefficient vector. We consider the high-dimensional setting, where
$\log(p)=o(n^b)$ for some constant $0<b<1$. We assume the distributions of $\bx$ and $\epsilon$ are
independent and both have mean 0. Under this assumption, $\betas$ is the mean effect of $y$
conditioning on $\bx$, which is assumed to be of interest.

To adapt for different magnitude of errors and robustify the estimation, we propose to use the
Huber loss \citep{HUB64}:
\begin{equation}
  \label{eq2.2}
  \ell_{\alpha}(x)=
  \begin{cases}
    2\alpha^{-1}|x|-\alpha^{-2}& \text{ if } |x|>\alpha^{-1}; \\
    x^2 & \text{ if } |x|\leq \alpha^{-1}.
  \end{cases}
\end{equation}
The Huber loss is quadratic for small values of $x$ and linear for large values of $x$. The parameter
$\alpha$ controls the blending of quadratic and linear penalization. The least squares and the LAD
can be regarded as two extremes of the Huber loss for $\alpha=0$ and $\alpha=\infty$,
respectively. Deviated from the traditional Huber's estimator, the parameter $\alpha$ converges to
zero in order to reduce the biases of estimating the mean regression function when the conditional
distribution of $\varepsilon_i$ is not symmetric.  On the other hand, $\alpha$ can not shrink too
fast in order to maintain the robustness. In this paper, we regard $\alpha$ as a tuning parameter,
whose optimal value will be discussed later in this section. In practice, $\alpha$ needs to be
tuned by some data-driven method. By letting $\alpha$ vary, we call $\ell_{\alpha}(x)$ the robust
approximate quadratic (RA-quadratic) loss.

To estimate $\betas$, we propose to solve the following convex optimization problem:
\begin{equation}
  \label{eq2.3}
  \hat{\bbeta}=\argmin_{\bbeta} \frac{1}{n} \sum_{i=1}^n
  \ell_{\alpha}(y_i-\bx_i^T\bbeta)+ \lambda_n \sum_{j=1}^p |\beta_j|.
\end{equation}
To assess the performance of $\hat{\bbeta}$, we study the property of
$\ltwonorm{\hat{\bbeta}-\betas}$, where $\ltwonorm{\cdot}$ is the Euclidean norm of a vector. When
$\lambda_n$ converges to zero sufficiently fast, $\hat{\bbeta}$ is a natural $M$-estimator of
$\betaa=\argmin_{\bbeta}\E \ell_{\alpha}(y-\bx'\bbeta)$, which is the population minimizer under
the RA-quadratic loss and varies by $\alpha$. In general, $\betaa$ differs from $\betas$.  But,
since the RA-quadratic loss approximates the quadratic loss as $\alpha$ tends to 0, $\betaa$ is
expected to converge to $\betas$. This property will be established in Theorem~\ref{thm1}.
Therefore, we decompose the statistical error $\hat{\bbeta}-\betas$ into the approximation error
$\betaa-\betas$ and the estimation error $\hat{\bbeta}-\betaa$. The statistical error
$\ltwonorm{\hat{\bbeta}-\betas}$ is then bounded by
\begin{equation*}
  \ltwonorm{\hat{\bbeta}-\betas}\leq
  \underbrace{\ltwonorm{\betaa-\betas}}_{\text{approximation error}} + \underbrace{\ltwonorm{\hat{\bbeta}-\betaa}}_{\text{estimation error}}.
\end{equation*}
In the following, we give the rate of the approximation and estimation error, respectively. We show
that $\ltwonorm{\hat{\bbeta}-\betas}$ admits the same rate as the optimal rate under light tails,
as long as the two tuning parameters $\alpha$ and $\lambda_n$ are  properly chosen. We first give
the rate of the approximation error under some moment conditions on $\bx$ and $\epsilon$. We assume
both $\betas$ and $\betaa$ are interior points of an $L_2$ ball with sufficiently large radius.

\begin{thm}
  \label{thm1}
  It holds that $\ltwonorm{\betaa-\betas}=O(\alpha^{k-1})$, under the
  following conditions:
  \begin{itemize} \itemsep -0.05in
  \item [(C1)] $\E|\epsilon|^k\leq M_k<\infty$, for some $k\geq 2$.
  \item [(C2)] $0<\kl\leq\lambda_{\min}(\E[\bx\bx^T])\leq\lambda_{\max}(\E[\bx\bx^T])\leq \ku<\infty$,
  \item [(C3)] For any $\bnu\in \Rp$, $\bx^T\bnu$ is sub-Gaussian with parameter
    at most $\ko^2 \ltwonorm{\bnu}^2$, i.e. $\E \exp (t\bx^T\bnu) \leq \exp(t^2\ko^2
    \ltwonorm{\bnu}^2/2)$, for any $t\in  \mathbb{R}$.
  \end{itemize}
\end{thm}

Theorem \ref{thm1} reveals that the approximation error vanishes faster if higher moments of error
distribution exists. We next give the rate of the estimation error
$\ltwonorm{\hat{\bbeta}-\betaa}$. This part differs from the existing work regarding the estimation
error of high dimensional regularized $M$-estimator \citep{NRW12,ANW12} as the population
minimizer $\betaa$ now varies with $\alpha$.
However, we will show that the estimation error rate does not depend
on $\alpha$, given a uniform sparsity condition.

In order to be solvable in the high-dimensional setting, $\betas$ is usually assumed to be sparse
or weakly sparse, i.e. many elements of $\betas$ are zero or small. By Theorem \ref{thm1}, $\betaa$
converges to $\betas$ as $\alpha$ goes to 0. In view of this fact, we assume that $\betaa$ is
uniformly weakly sparse when $\alpha$ is sufficiently small. In particular, we assume that there
exists a small constant $r>0$, such that $\betaa$ belongs to an $L_q$-ball with a uniform radius
$\Rq$ that
\begin{equation}
  \label{eq2.4}
  \sum_{j=1}^p |\beta_{\alpha,j}^{\ast}|^q\leq \Rq,
\end{equation}
for all $\alpha\in (0,r]$, and some $q\in [0,1]$.  When the conditional distribution of $\varepsilon_i$ is symmetric, $\beta_{\alpha,j}^{\ast} = \beta_{j}^*$ for all $\alpha$ and $j$.  Therefore the condition reduces to that $\betas$ is in the $L_q$ ball. In a special case where $q=1$, it follows from
Theorem \ref{thm1} that if $\betas$ belongs to the $L_1$-ball with radius $R_1/2$ and $r\leq
[R_1/(2c_0 \sqrt{p})]^{\frac{1}{k-1}}$, where $c_0$ is a generic constant, then $\betaa$ belongs to
the $L_1$-ball with radius $R_1$ for all $\alpha\in (0,r]$. For a general $q\in[0,1)$, we assume a
uniform upper bound $\Rq$ as in (\ref{eq2.4}), which is allowed to diverge to infinity.

Since the RA-quadratic loss is convex, we show that with high probability the estimation error
$\hat{\bDelta}=\hat{\bbeta}-\betaa$ belongs to a star-shaped set, which depends on $\alpha$ and the
threshold level $\eta$ of signals.

\begin{lem}
  \label{lem1} Under Conditions (C1) and (C3), by
  choosing $\lambda_n=\kappa_{\lambda}\sqrt{\frac{\log p}{n}}$ and $\alpha\geq
  \frac{L\lambda_n}{4v}$, where $\kappa_{\lambda}$, $v$ and $L$ are some constants, with
  probability greater than $1-2p^{-c_0}$,
  \begin{equation*}
    \hat{\bDelta}=\hat{\bbeta}-\betaa\in \bbC_{\alpha \eta}=\{\bDelta\in\Rp:
    \lonenorm{\bDelta_{S_{\alpha\eta}^c}}\leq 3 \lonenorm{\bDelta_{S_{\alpha\eta}}}+4
    \lonenorm{\bbeta_{S_{\alpha \eta}^c}^{\ast}} \},
  \end{equation*}
  where $c_0=\kappa_{\lambda}^2/(32 v)-1$, $\eta$ is a positive constant,
  $S_{\alpha\eta}=\{j:|\beta_{\alpha,j}^{\ast}|>\eta \}$ and $\bDelta_{S_{\alpha \eta}}$ denotes
  the subvector of $\bDelta$ with indices in set $S_{\alpha \eta}$.
\end{lem}

We further verify a restricted strong convexity (RSC) condition, which has been shown to be
critical in the study of high dimensional regularized $M$-estimator
\citep{NRW12,ANW12}. Let
\begin{equation}
  \label{eq2.5}
  \delta\cL_n(\bDelta,\bbeta)=\cL_n(\bbeta+\bDelta)-\cL_n(\bbeta)-[\nabla \cL_n(\bbeta)]^T\bDelta,
\end{equation}
where $\cL_n(\bbeta)=\frac{1}{n}\sum_{i=1}^n \ell_{\alpha}(y_i-\bx_i^T\bbeta)$, $\bDelta$ is a
$p$-dimensional vector and $\nabla\cL_n(\bbeta)$ is the gradient of $\cL_n$ at the point of $\bbeta$.
\begin{dfn}
  \label{def:RSC}
  The loss function $\cL_n$ satisfies RSC condition on a set $S$ with curvature $\kappa_{\cL}>0$ and
  tolerance $\tau_{\cL}$ if
  \begin{equation*}
    \delta \cL_n(\bDelta, \bbeta)\geq \kappa_{\cL}\ltwonorm{\bDelta}^2-\tau_{\cL}^2, \text{ for all
    } \bDelta\in S.
  \end{equation*}
\end{dfn}
Next, we show that with high probability the RA-quadratic loss (\ref{eq2.2}) satisfies RSC for all
$\bDelta \in \bbC_{\alpha \eta} \cap \{\bDelta:\ltwonorm{\bDelta}\leq 1 \}$ with uniform constants $\kappa_{\cL}$ and $\tau_{\cL}$ that do not
depend on $\alpha$. Lemma \ref{lem2} is a preliminary result, based on which RSC is checked in
Lemma \ref{lem3}.

\begin{lem}
  \label{lem2}
  Under conditions (C1)-(C3), for all $\ltwonorm{\bDelta}\leq 1$, there exist uniform positive
  constants $\kappa_1$ and $\kappa_2$, such that
  \begin{equation}
    \label{eq2.6}
    \delta \cL_n(\bDelta,\betaa)\geq \kappa_1 \ltwonorm{\bDelta}(\ltwonorm{\bDelta}-\kappa_2 \sqLogpN
    \lonenorm{\bDelta}),
  \end{equation}
  with probability at least $1-c_1\exp(-c_2n)$ for some positive constants $c_1$ and $c_2$.
\end{lem}

\begin{lem}
  \label{lem3}
  Suppose conditions  (C1)-(C3) hold and assume that
  \begin{equation}
    \label{eq2.7}
    8\kappa_2\kappa_{\lambda}^{-q/2} \sqrt{\Rq} \left(\frac{\log p}{n}
    \right)^{(1-q)/2} \leq 1,
  \end{equation}
  with the choice $\eta=\lambda_n$, with probability at least $1-c_1\exp(-c_2n)$, the RSC condition
  holds for $\delta\cL_n(\bDelta,\betaa)$ and $\bDelta\in\bbC_{\alpha \eta}\cap
  \{\bDelta:\ltwonorm{\bDelta}\leq 1 \}$ with
  $\kappa_{\mathcal{L}}={\kappa_1}/{2}$ and $\tau_{\mathcal{L}}^2=4\Rq\kappa_2
  \kappa_{\lambda}^{1-q} \left(\frac{\log p}{n} \right)^{1-(q/2)}$.
\end{lem}

Lemma \ref{lem3} shows that, even though $\betaa$ is unknown and the set $\bbC_{\alpha\eta}$
depends on $\alpha$, RSC holds with uniform constants that do not depend on $\alpha$. This further
gives the following upper bound of the estimation error $\ltwonorm{\hat{\bbeta}-\betaa}$, which
also does not depend on $\alpha$.

\begin{thm}
  \label{thm2}
  Under conditions of Lemma \ref{lem1} and \eqref{eq2.7}, with probability at least
  $1-2p^{-c_0}-c_1\exp(-c_2n)$,
  \begin{equation*}
    \ltwonorm{\hat{\bbeta}-\betaa}=O(\sqrt{\Rq}[(\log p)/n]^{1/2-q/4}).
  \end{equation*}
\end{thm}
Finally, Theorem \ref{thm1} and \ref{thm2} together lead to the following main result, which
gives the rate of the statistical error $\ltwonorm{\hat{\bbeta}-\betas}$.
\begin{thm}
  \label{thm3}
  Under conditions of Lemma \ref{lem1} and \eqref{eq2.7}, with probability at least
  $1-2p^{-c_0}-c_1\exp(-c_2n)$,
  \begin{equation*}
    \ltwonorm{\hat{\bbeta}-\betas}=O({\alpha}^{k-1})+O(\sqrt{\Rq}[(\log p)/n]^{1/2-q/4}).
  \end{equation*}
\end{thm}

Next, we compare our result with the existing results regarding the robust estimation of high
dimensional linear regression model.
\begin{enumerate}
\item When the  distribution of $\epsilon$ is symmetric around 0, then
  $\betaa=\betas$ for any $\alpha$, which has no approximation error. If
  $\epsilon$ has heavy tails in addition to being symmetric, we
  would like to choose $\alpha$ sufficiently large to robustify the estimation. It then follows
  from Theorem \ref{thm2} that $\ltwonorm{\hat{\bbeta}-\betas}=O_P(\sqrt{R_q}[(\log
  p)/n]^{1/2-q/4})$, where $R_q=\sum_{j=1}^p |\beta_j^{\ast}|^q$. The rate is the same as the
  minimax rate \citep{RWY11} for weakly sparse model under the light tails. In a special case that
  $q=0$, it gives $\ltwonorm{\hat{\bbeta}-\betas}=O_P(\sqrt{s(\log p)/n})$, where $s$ is the number
  of nonzero elements in $\betas$. This is the same rate as the regularized LAD estimator in
  \cite{WAN13} and the regularized quantile estimator in \cite{BC11}. It suggests that our method
  does not lose efficiency for symmetric heavy-tailed errors.

\item If the distribution of $\epsilon$  is asymmetric around 0, the quantile and
  LAD based methods are inconsistent, since they estimate the median instead of the mean. Theorem \ref{thm3}
  shows that our estimator still achieves the optimal rate given that $\alpha=o(\{R_q[(\log
  p)/n]^{1-\frac{q}{2}}\}^{\frac{1}{2(k-1)}})$ even though the $k$-th moment of $\varepsilon$ is assumed. Recall from conditions in Lemma \ref{lem1} that we
  also need to choose $\alpha>c \sqrt{(\log p)/n}$ for some constant $c$. Given the sparsity
  condition (\ref{eq2.7}), $\alpha$ can be chosen to meet the above two requirements. In terms of
  estimating the conditional mean effect, errors with heavy but asymmetric tails give the case
  where the RA-Lasso has the biggest advantage over the other estimators.
\end{enumerate}
In practice, the distribution of errors is unknown. However, we proved that our method is no worse
than the existing methods for any type of errors, as long as the tuning parameters are chosen
properly. Hence, our method is more flexible.

\section{Geometric convergence of computational error}
\label{sec3}
The gradient descent algorithm \citep{N07,ANW12} is usually applied to solve the convex problem
(\ref{eq2.3}). For example, we can replace the RA-quadratic loss with its local quadratic
approximation (LQA) and iteratively solve the following optimization problem:
\begin{equation}
  \label{eq3.1}
  \hat{\bbeta}^{t+1}=\argmin_{\lonenorm{\bbeta}\leq\rho}
  \left\{\cL_n(\hat{\bbeta}^t)+[\nabla\cL_n(\hat{\bbeta}^t)]^T (\bbeta-\hat{\bbeta}^t)
    +\frac{\gamma_{u}}{2} \ltwonorm{\bbeta-\hat{\bbeta}^t}^2+\lambda_n \lonenorm{\bbeta}\right\},
\end{equation}
where $\gamma_{u}$ is a fixed constant at each iteration, and the side constraint
``$\lonenorm{\bbeta}\leq \rho$'' is introduced to guarantee good performance in the first few
iterations and $\rho$ is allowed to be sufficiently large. To solve (\ref{eq3.1}), the update can
be computed by a two-step procedure. We first solve (\ref{eq3.1}) without the norm constraint by
soft-thresholding the vector $\hat{\bbeta}^t - \frac{1}{\gamma_{u}}\nabla\cL_n(\hat{\bbeta}^t)$ at
level $\lambda_n$ and call the solution $\check{\bbeta}$. If $\lonenorm{\check{\bbeta}}\leq \rho$,
set $\hat{\bbeta}^{t+1}=\check{\bbeta}$. Otherwise, $\hat{\bbeta}^{t+1}$ is obtained by further
project $\check{\bbeta}$ onto the $L_1$-ball $\{\bbeta: \lonenorm{\bbeta} \leq \rho\}$. The
projection can be done \citep{DSS08} by soft-thresholding $\check{\bbeta}$ at level $\pi_n$, where
$\pi_n$ is given by the following procedure: (1) sort $\{|\check{\beta}_j| \}_{j=1}^p$ into
$b_1\geq b_2\geq \ldots\geq b_p$; (2) find $J=\max\{1\leq j\leq p: b_j-
(\sum_{r=1}^{j}b_r-\rho)/j>0\}$ and let $\pi_n=(\sum_{r=1}^J b_j-\rho)/J$.

\cite{ANW12} considered the computational error of such first-order gradient descent method. They
showed that, for a convex and differentiable loss functions $\ell(x)$ and decomposable penalty
function $p(\bbeta)$, the error $\ltwonorm{\hat{\bbeta}^t-\betas}$ has the same rate as
$\ltwonorm{\hat{\bbeta}-\betas}$ for all sufficiently large $t$, where $\betas=\argmin_{\bbeta}
\E\cL(\bx,y;\bbeta)$, and $\hat{\bbeta}=\argmin_{\bbeta} \frac{1}{n} \sum_{i=1}^n
\cL(\bx_i,y_i,\bbeta)+p(\bbeta)$. Different from their setup, our population minimizer $\betaa$
varies by $\alpha$. Nevertheless, as $\betaa$ converges to the true effect $\betas$, by a careful
control of $\alpha$, we can still show that $\ltwonorm{\hat{\bbeta}^t-\betas}$ has the same rate as
$\ltwonorm{\hat{\bbeta}-\betas}$, where $\hat{\bbeta}$ is the theoretical solution of (\ref{eq2.3})
and $\hat{\bbeta}^t$ is as defined in (\ref{eq3.1}).

The key is that the RA-quadratic loss function $\cL_n$ satisfies the restricted strong convexity (RSC)
condition and the restricted smoothness condition (RSM) with some uniform constants, namely
$\delta\cL_n(\bDelta,\bbeta)$ as defined in (\ref{eq2.5}) satisfies the following conditions:
\begin{align}
  &\text{RSC}:\delta \cL_n(\bDelta,\bbeta)\geq
  \frac{\gamma_{l}}{2}\ltwonorm{\bDelta}^2-\tau_{l}\lonenorm{\bDelta}^2, \label{eq3.2}\\
  &\text{RSM}:\delta \cL_n(\bDelta,\bbeta)\leq
  \frac{\gamma_{u}}{2}\ltwonorm{\bDelta}^2+\tau_{u}\lonenorm{\bDelta}^2, \label{eq3.3}
\end{align}
for all $\bbeta$ and $\bDelta$ in some set of interest, with parameters $\gamma_l$, $\tau_l$,
$\gamma_u$ and $\tau_u$ that do not depend on $\alpha$. We show that such conditions hold with high
probability.
\begin{lem}
  \label{lem4}
  Under conditions (C1)-(C3), for all $\bbeta\in \mathbb{R}^p$ and $\bDelta \in
  \{\bDelta:\ltwonorm{\bDelta}\leq 1 \}$, with probability greater
  than $1-c_1\exp(-c_2n)$, (\ref{eq3.2}) and (\ref{eq3.3}) hold with $\gamma_l=\kappa_1$,
  $\tau_l=\kappa_1\kappa_2^2(\log p)/(2n)$, $\gamma_u=3\ku$, $\tau_u=\ku(\log p)/n$.
\end{lem}

We further show in Theorem \ref{thm4} that, whenever $\Rq(\frac{\log p}{n})^{1-(q/2)}=o(1)$, which
is required for consistency of \emph{any method} over the weak sparse $L_q$ ball by the known
minimax results \citep{RWY11}, it holds that
$\ltwonorm{\hat{\bbeta}^t-\hat{\bbeta}}=o(\ltwonorm{\hat{\bbeta}-\betaa})$ for sufficiently many
iterations with lower bound specified in Theorem~\ref{thm4}. Hence,
\begin{align*}
  \ltwonorm{\hat{\bbeta}^t-\betas} & \leq \ltwonorm{\hat{\bbeta}^t-\hat{\bbeta}}+\ltwonorm{\hat{\bbeta}
    -\betaa}+\ltwonorm{\betaa-\betas}\\
  & =o(\ltwonorm{\hat{\bbeta}
    -\betaa})+\ltwonorm{\hat{\bbeta}-\betaa}+\ltwonorm{\betaa-\betas}\\
  & =O({\alpha}^{k-1})
  +O(\sqrt{\Rq}[(\log p)/n]^{1/2-q/4}),
\end{align*}
which has the same rate as
$\ltwonorm{\hat{\bbeta}-\betas}$. Hence, from a statistical point of view, there is no need to
iterate beyond $t$ steps.

\begin{thm}
  \label{thm4}
  Under conditions of Theorem \ref{thm3}, suppose we choose $\lambda_n$ as in Lemma \ref{lem1}
  and also satisfying
  \begin{equation*}
    \lambda_n\geq \frac{32\rho}{1-\kappa}\left(1-\frac{64\ku|S_{\alpha\eta}|\frac{\log
          p}{n}}{\bar{\gamma}_l}\right)^{-1} \left[1+\kappa_1\kappa_2^2
      \left(\frac{\bar{\gamma}_l}{12\ku}+\frac{128\ku|S_{\alpha\eta}|\frac{\log
            p}{n}}{\bar{\gamma}_l} \right)+8\ku\right]\frac{\log p}{n},
  \end{equation*}
  where $|S_{\alpha\eta}|$ denotes the cardinality of set $S_{\alpha\eta}$ and
  $\bar{\gamma}_l=\gamma_l-64\tau_l|S_{\alpha \eta}|$, then with probability at least
  $1-p^{-c_0}-c_1\exp(-c_2n)$, we have
  \begin{equation*}
    \ltwonorm{\hat{\bbeta}^t-\hat{\bbeta}}^2=O \left(\Rq \left(\frac{\log
          p}{n}\right)^{1-(q/2)}\left[\ltwonorm{\hat{\bbeta}-\betaa}^2+\Rq \left(\frac{\log
            p}{n}\right)^{1-(q/2)}\right] \right),
  \end{equation*}
  for all iterations
  \begin{equation*}
    t\geq \frac{2\log((\phi_n(\hat{\bbeta}^0)-\phi_n(\hat{\bbeta}))/\delta^2)}{\log(1/\kappa)}+
    \log_2\log_2 \left(\frac{\rho\lambda_n}{\delta^2}
    \right)\left(1+\frac{\log2}{\log(1/\kappa)}\right),
  \end{equation*}
  where $\phi_n(\bbeta)=\cL_n(\bbeta)+\lambda_n \lonenorm{\bbeta}$ and $\hat{\bbeta}^0$ is the
  initial value, $\delta=\varepsilon^2/(1-\kappa)$ is the tolerance level, $\kappa$ and
  $\varepsilon$ are some constants as will be defined in (\ref{eq7.21}) and
  (\ref{eq7.22}), respectively.
\end{thm}

\section{Connection with Catoni loss}
\label{sec4}
\cite{CAT12} considered the estimation of the mean of heavy-tailed distribution with fast
concentration. He proposed an $M$-estimator by solving
\begin{equation*}
  \sum_{i=1}^n \psi_c[\alpha(y_i-\theta)]=0,
\end{equation*}
where the influence function $\psi_c(x)$ is chosen such that $-\log(1-x+x^2/2)\leq \psi_c(x)\leq
\log(1+x+x^2/2)$. He showed that this $M$-estimator has the exponential type of concentration by
only requiring the existence of the variance. It performed as well as the sample mean under the light-tail
case. In Section \ref{sec2}, we essentially showed the same type of
concentration for the RA-quadratic loss under the linear regression setting.

The estimation of mean can be regarded as a univariate linear regression where the covariate equals
to 1. In that special case, we have a more explicit concentration result for the RA-mean estimator,
which is the estimator that minimizes the RA-quadratic loss. Let $\{y_i \}_{i=1}^n$ be an i.i.d
sample from some unknown distribution with $\E(y_i)=\mu$ and $\var(y_i)=\sigma^2$. The RA-mean
estimator $\hat{\mu}_{\alpha}$ of $\mu$ is the solution of
\begin{equation*}
  \sum_{i=1}^n \psi[\alpha(y_i-\mu)]=0,
\end{equation*}
for parameter $\alpha \to 0$, where the influence function $\psi(x)=x$ if $|x|\leq 1$, $\psi(x)=1$,
if $x>1$ and $\psi(x)=-1$ if $x<-1$. The following theorem gives the exponential type of
concentration of $\hat{\mu}_{\alpha}$ around $\mu$.

\begin{thm}
  \label{thm5}
  Assume $\frac{\log (1/\delta)}{n}\leq 1/8$ and let $\alpha=\sqrt{\frac{\log
      (1/\delta)}{nv^2}}$ where $v \geq \sigma$.  Then,
  \begin{equation*}
    P\left ( |\hat{\mu}_{\alpha}-\mu|\geq 4v \sqrt{\frac{\log (1/\delta)}{n}} \right ) \leq 2 \delta.
  \end{equation*}
\end{thm}

The above result provides fast
concentration of the mean estimation with only two moments assumption.
This is very useful for last scale hypothesis testing \citep{EFR10, FHG12}
and covariance matrix estimation \citep{BL08, FLM13}, where uniform convergence is required.
Taking the estimation of large covariance matrix as an example, in order for the elements of the
sample covariance matrix to converge uniformly, the aforementioned authors require the underlying
multivariate distribution be sub-Gaussian.  This restrictive assumptions can be removed if we apply
the robust estimation with concentration bound.  Regarding $\sigma_{ij} = \E X_i X_j$ as the
expected value of the random variable $X_i X_j$ (it is typically not the same as the median of $X_i
X_j$), it can be estimated with accuracy
$$
P\left ( |\hat{\sigma}_{ij}-\sigma_{ij}|\geq 4v \sqrt{\frac{\log (1/\delta)}{n}} \right ) \leq 2
\delta,
$$
where $v \geq \max_{i,j\leq  p} \sqrt{\var(X_iX_j)}$ and $\hat{\sigma}_{ij}$  is RA-mean estimator
using data $\{X_{ik} X_{jk}\}_{k=1}^n$.  Since there  are only $O(p^2)$ elements, by taking $\delta
= p^{-3}$ and the union bound, we have
$$
\max_{i, j \leq p} \sqrt{n/\log p} |\hat{\sigma}_{ij}-\sigma_{ij}| \to 0,
$$
when $\E X_i^4<\infty$. This robustified covariance estimator requires much weaker condition
than the sample covariance and has far wide applicability than the sample covariance.  It can be
regularized further in the same way as the sample covariance matrix.

On the other hand, Catoni's idea could also be extended to the linear regression setting. Suppose
we replace the RA-quadratic loss $\ell_{\alpha}(x)$ in (\ref{eq2.3}) with Catoni loss
\begin{equation*}
  \ell^c_{\alpha}(x)=\frac{2}{\alpha}\int_0^x\psi_c(\alpha t)dt,
\end{equation*}
where the influence function $\psi_c(t)$ is given by
\begin{equation*}
  \psi_c(t)= \mbox{sgn}(t) \{-\log(1-|t|+t^2/2) I(|t| < 1) + \log(2) I(|t| \geq 1) \}.
\end{equation*}
Let $\hat{\bbeta}^c$ be the corresponding solution. Then, $\hat{\bbeta}^c$ has the same convergence
rate as the RA-Lasso, when the second or the third moment of errors exists.

\begin{thm}
  \label{thm6}
  Suppose condition (C1) holds for $k=2$ or 3, (C2), (C3) and (\ref{eq2.7}) hold, then with probability at least
  $1-2p^{-c_0}-c_1\exp(-c_2n)$,
  \begin{equation*}
    \ltwonorm{\hat{\bbeta}^c-\betas}=O({\alpha}^{k-1})+O(\sqrt{\Rq}[(\log p)/n]^{1/2-q/4}).
  \end{equation*}
\end{thm}

Unlike the RA-lasso, the order of bias of $\hat{\bbeta}^c$ cannot be further improved, even when higher moments of errors exist
beyond the third order. The reason is that the Catoni loss is not exactly the quadratic loss over any finite
intervals. Similar results regarding the computational error of $\hat{\bbeta}^c$ could also be established as in
Theorem \ref{thm4}, since the RSC/RSM conditions also hold for Catoni loss with uniform constants.

\section{Variance Estimation}
\label{sec5}
We estimate $\sigma^2=E\epsilon^2$ based on the RA-Lasso estimator and a cross-validation scheme. To ease the
presentation, we assume the data set can be evenly divided into $K$ folds with $m$ observations in each
fold. Then, we estimate $\sigma^2$ by
\begin{equation*}
  \hat{\sigma}^2=\frac{1}{K} \sum_{k=1}^K \frac{1}{m} \sum_{i\in \text{fold }k}
  (y_i-\bx_i^T\hat{\bbeta}^{(-k)})^2,
\end{equation*}
where $\hat{\bbeta}^{(-k)}$ is the RA-Lasso estimator obtained by using data points outside the
$k$-th fold. We show that $\hat{\sigma}^2$ is asymptotically efficient.
\begin{thm}
  \label{thm7}
  Under conditions of Theorem \ref{thm3}, if $\Rq (\log p)^{1-q/2}/n^{(1-q)/2} \to 0$ for $q\in
  [0,1)$, and  $\alpha=o\left(\{R_q[(\log
    p)/n]^{1-\frac{q}{2}}\}^{\frac{1}{2(k-1)}}\right)$, then
  \[\sqrt{n}(\hat{\sigma}^2 - \sigma^2)\xrightarrow{\mathcal{D}} {N}(0, \E\epsilon^4 - \sigma^4).\]
\end{thm}

\section{Simulation Studies}
\label{sec6}
In this section, we assess the finite sample performance of the RA-Lasso and compare it with other
methods through various models. We simulated data from the following high dimensional model
\begin{equation} \label{eq6.1}
  y_i = \bx_i^T \betas + \epsilon_i, \quad \bx_i \sim N(0, I_p),
\end{equation}
where we generated $n=100$ observations and the number of parameters was chosen to be $p=400$. We
chose the true regression coefficient vector as
\begin{equation*}
  \bbeta^{\ast} = (3, \ldots, 3, 0, \ldots , 0)^T,
\end{equation*}
where the first 20 elements are all equal to 3 and the rest are all equal to 0.  To involve various
shapes of error distributions, we considered the following five scenarios:
\begin{enumerate}
\item Normal with mean 0 and variance 4 (N(0,4));
\item Two times the t-distribution with degrees of freedom 3 ($2t_3$);
\item Mixture of Normal distribution(MixN): $0.5N(-1,4)+0.5N(8,1)$;
\item Log-normal distribution (LogNormal): $\epsilon =e^{1+1.2Z}$, where $Z$ is standard normal.
\item Weibull distribution with shape parameter = 0.3 and scale parameter = 0.5 (Weibull).
\end{enumerate}
In order to meet the model assumption, the errors were standardized to have mean 0. Table \ref{tab1} categorizes the five scenarios according to the shapes and tails of the error distributions.
\begin{table}[htp!]
  \begin{center}
    \begin{tabular}{c|cc}
      \hline\hline
      & \textbf{Light Tail} & \textbf{Heavy Tail}\\
      \hline
      \textbf{Symmetric} & $N(0,4)$ & $2t_3$ \\
      \textbf{Asymmetric} & MixN & LogNormal, Weibull \\
      \hline\hline
    \end{tabular}
  \end{center}
  \caption{Summary of the shapes and tails of five error distributions \label{tab1}}
\end{table}

To obtain our estimator, we iteratively applied the gradient descent algorithm. We compared
RA-Lasso with another two methods in high-dimensional setting: (a) Lasso: the penalized
least-squares estimator with $L_1$-penalty as in \cite{TIB96}; and (b) R-Lasso: the R-Lasso
estimator in \cite{FFB14}, which is the same as the regularized LAD estimator with $L_1$-penalty as
in \cite{WAN13}. Their performance under the five scenarios was evaluated by the following four
measurements:
\begin{itemize}\itemsep -0.05 in
\item [(1)] $L_2$ error, which is defined as $\ltwonorm{\hat{\bbeta}-\betas}$.
\item [(2)] $L_1$ error, which is defined as $\lonenorm{\hat{\bbeta}-\betas}$.
\item [(3)] Number of false positives (FP), which is number of noise covariates that are selected.
\item [(4)] Number of false negatives (FN), which is number of signal covariates that are not selected.
\end{itemize}
We also measured the relative gain of RA-Lasso with respect to R-Lasso and Lasso, in terms of the
difference to the oracle estimator. The oracle estimator $\hat{\bbeta}_{\text{oracle}}$ is defined
to be the least square estimator by using the first 20 covariates only. Then, the relative gain of
RA-Lasso with respect to Lasso ($\text{RG}_{\text{A,L}}$) in $L_2$ and $L_1$ norm are defined as
\begin{equation*}
  \frac{\ltwonorm{\hat{\bbeta}_{\text{Lasso}}-\betas}-\ltwonorm{\hat{\bbeta}_{\text{oracle}}-\betas}}{\ltwonorm{\hat{\bbeta}_{\text{RA-Lasso}}-\betas}-\ltwonorm{\hat{\bbeta}_{\text{oracle}}-\betas}}
  \quad
  \text{ and }
  \quad
  \frac{\lonenorm{\hat{\bbeta}_{\text{Lasso}}-\betas}-\lonenorm{\hat{\bbeta}_{\text{oracle}}-\betas}}{\lonenorm{\hat{\bbeta}_{\text{RA-Lasso}}-\betas}-\lonenorm{\hat{\bbeta}_{\text{oracle}}-\betas}}.
\end{equation*}
The relative gain of RA-Lasso with respect to R-Lasso ($\text{RG}_{\text{A,R}}$) is defined
similarly.

For RA-Lasso, the tuning parameters $\lambda_n$ and $\alpha$ were chosen optimally based on 100
independent validation datasets. We ran a 2-dimensional grid search to find the best
$(\lambda_n,\alpha)$ pair that minimizes the mean $L_2$-loss of the 100 validation datasets. Such
an optimal pair was then used in the simulations. Similar method was applied in choosing the tuning
parameters in Lasso and R-Lasso.

The above simulation model is based on the additive model \eqref{eq6.1}, in which error
distribution is independent of covariates. However, this homoscedastic model makes the conditional
mean and the conditional median differ only by a constant. To further examine the deviations between
the mean regression and median regression, we also simulated the data from the heteroscedastic
model
\begin{equation} \label{eq6.2}
  y_i = \bx_i^T \betas + c^{-1} (\bx_i^T \betas)^2 \epsilon_i, \quad \bx_i \sim N(0, I_p),
\end{equation}
where the constant $c = \sqrt{3}\| \betas \|^2$ makes $ \E [c^{-1} (\bx_i^T \betas)^2]^2 = 1$.  Note that $\bx_i^T
\betas \sim N(0, \|\betas\|^2)$ and therefore $c$ is chosen so that the average noise levels is the same as that
of $\epsilon_i$. For both the homoscedastic and the heteroscedastic models, we ran 100 simulations for each
scenario. The mean of each performance measurement is reported in Table \ref{tab2} and \ref{tab3}, respectively.
\begin{table}[hbtp]
  \begin{center}
    \begin{tabular}{cc|ccccc}
      \hline\hline
      &            & \textbf{Lasso} & \textbf{R-Lasso} & \textbf{RA-Lasso} & \textbf{$\text{RG}_{\text{A,L}}$} & \textbf{$\text{RG}_{\text{A,R}}$}\\
      \hline
      \multirow{3}{*}{$\mathbf{N(0,4)}$}  & $L_2$ loss & 4.54           & 4.40             & 4.53              & 1.00                              & 0.96 \\

      & $L_1$ loss & 27.21          & 29.11            & 27.21             & 1.00                              & 1.08\\

      & FP, FN     & 52.10, 0.09    & 66.36, 0.17      & 52.10, 0.09       &                                   & \\
      \hline
      \multirow{3}{*}{$\mathbf{2t_3}$}    & $L_2$ loss & 6.04           & 5.10             & 5.47              & 1.14                              & 0.91 \\

      & $L_1$ loss & 35.22          & 33.07            & 30.42             & 1.19                              & 1.10\\

      & FP, FN     & 47.13, 0.34    & 65.84, 0.22      & 41.34, 0.28       &                                   & \\

      \hline
      \multirow{3}{*}{\textbf{MixN}}      & $L_2$ loss & 6.14           & 6.44             & 6.13              & 1.00                              & 1.06\\
      & $L_1$ loss & 40.46          & 46.18            & 38.48             & 1.06                              & 1.23\\
      & FP, FN     & 65.99, 0.34    & 80.31, 0.33      & 58.05, 0.39       &                                   & \\

      \hline
      \multirow{3}{*}{\textbf{LogNormal}} & $L_2$ loss & 11.08          & 12.16            & 10.10             & 1.14                              & 1.30\\
      & $L_1$ loss & 53.17          & 57.18            & 51.58             & 1.04                              & 1.14\\
      & FP, FN     & 26.5, 15.00    & 27.20, 6.90      & 37.20, 3.90       &                                   & \\
      \hline
      \multirow{3}{*}{\textbf{Weibull}}   & $L_2$ loss & 7.77           & 7.11             & 6.62              & 1.23                              & 1.10\\
      & $L_1$ loss & 55.65          & 50.49            & 42.93             & 1.34                              & 1.20\\
      & FP, FN     & 78.70, 0.71    & 77.13, 0.56      &
      62.27,0.52                          &            & \\
      \hline\hline
    \end{tabular}
    \caption{Simulation results of Lasso, R-Lasso and RA-Lasso under homoscedastic model
      (\ref{eq6.1}).  \label{tab2}}
  \end{center}
\end{table}

\begin{table}[hbtp]
  \begin{center}
    \begin{tabular}{cc|ccccc}
      \hline\hline
      &            & \textbf{Lasso} & \textbf{R-Lasso} & \textbf{RA-Lasso} & \textbf{$\text{RG}_{\text{A,L}}$} & \textbf{$\text{RG}_{\text{A,R}}$}\\
      \hline
      \multirow{3}{*}{$\mathbf{N(0,4)}$}  & $L_2$ loss & 4.60           & 4.34             & 4.60              & 1.00              & 0.93\\
      & $L_1$ loss & 27.16          & 27.14            & 27.15             & 1.00              & 1.00 \\
      & FP, FN     & 48.78, 0.10    & 58.25, 0.27      & 48.78, 0.10       & \\
      \hline
      \multirow{3}{*}{$\mathbf{2t_3}$}    & $L_2$ loss & 8.08           & 6.71             & 6.70              & 1.26              & 1.01\\
      & $L_1$ loss & 41.16          & 42.76            & 38.52             & 1.08              & 1.12\\
      & FP, FN     & 55.33, 0.67    & 71.67, 0.33      & 45.33, 0.33       &                   & \\
      \hline
      \multirow{3}{*}{\textbf{MixN}}      & $L_2$ loss & 6.26           & 6.54             & 6.25              & 1.00              & 1.06\\
      & $L_1$ loss & 41.26          & 46.95            & 39.25             & 1.06              & 1.23\\
      & FP, FN     & 65.98, 0.34    & 80.30, 0.32       & 58.80 0.34        &                   & \\
      \hline
      \multirow{3}{*}{\textbf{LogNormal}} & $L_2$ loss & 10.86          & 9.19             & 8.48              & 1.43              & 1.13\\
      & $L_1$ loss & 57.52          & 57.18            & 53.20             & 1.10              & 1.09\\
      & FP, FN     & 29.70, 5.70      & 54.10, 2.00        & 54.30, 1.50         &                   & \\
      \hline
      \multirow{3}{*}{\textbf{Weibull}}   & $L_2$ loss & 7.40           & 8.81             & 5.53              & 1.53              & 1.92\\
      & $L_1$ loss & 40.95          & 47.82            & 34.65             & 1.23              & 1.48\\
      & FP, FN     & 38.87,0.96     & 35.31, 2.90      & 58.15,0.39        &                   & \\
      \hline\hline
    \end{tabular}
    \caption{Simulation results of Lasso, R-Lasso and RA-Lasso under heteroscedastic model
      (\ref{eq6.2}).  \label{tab3}}
  \end{center}
\end{table}

Tables \ref{tab2} and \ref{tab3} indicate that our method had the biggest advantage when the errors
were asymmetric and heavy-tailed (LogNormal and Weibull). In this case, R-Lasso had larger $L_1$
and $L_2$ errors due to the bias for estimating the conditional median instead of the mean. Even
though Lasso did not have bias in the loss component, it did not perform well due to its sensitivity to outliers. The
advantage of our method is more pronounced in the heteroscedastic model than in the homoscedastic
model. Both of them clearly indicate that if the errors come from asymmetric and heavy-tailed
distributions, our method is better than both Lasso and R-Lasso. When the errors were symmetric and
heavy-tailed ($2t_3$), our estimator performed closely as R-Lasso, both of which outperformed
Lasso. The above two cases evidently showed that RA-Lasso was robust to the outliers and did not
lose efficiency when the errors were indeed symmetric. Under the light-tailed scenario, if the
errors were asymmetric (MixN), our method performed similarly as Lasso. R-Lasso performed worse,
since it had bias. For the regular setting (N(0, 4)), where the errors were light-tailed and
symmetric, the three methods were comparable with each other.

In conclusion, RA-Lasso is more flexible than Lasso and R-Lasso. The tuning parameter $\alpha$
automatically adapts to errors with different shapes and tails. It enables RA-Lasso to render
consistently satisfactory results under all scenarios.

\section{Real Data Example}
\label{sec7}
\begin{figure}[hbtp!]
  \centering
  \includegraphics[width=\textwidth]{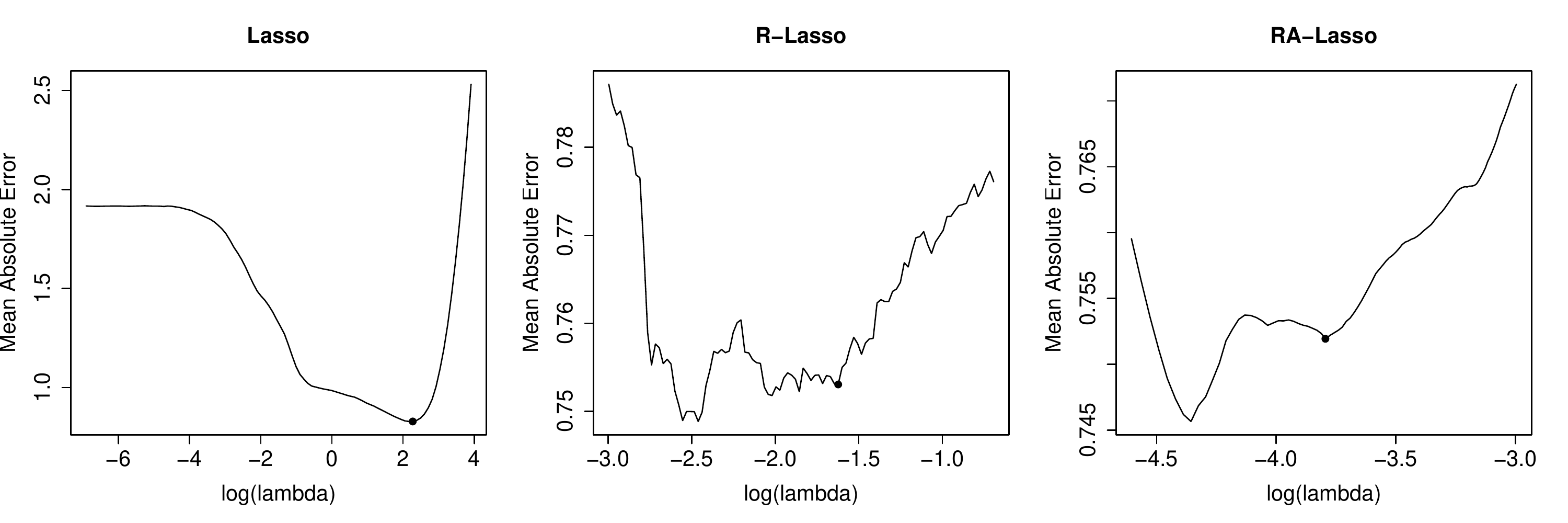}
  \caption{Five-fold cross validation results: black dot marks the choice of the penalization
    parameter. \label{fig1} }
\end{figure}

\begin{figure}[hbtp!]
  \centering
  \includegraphics[width=\textwidth]{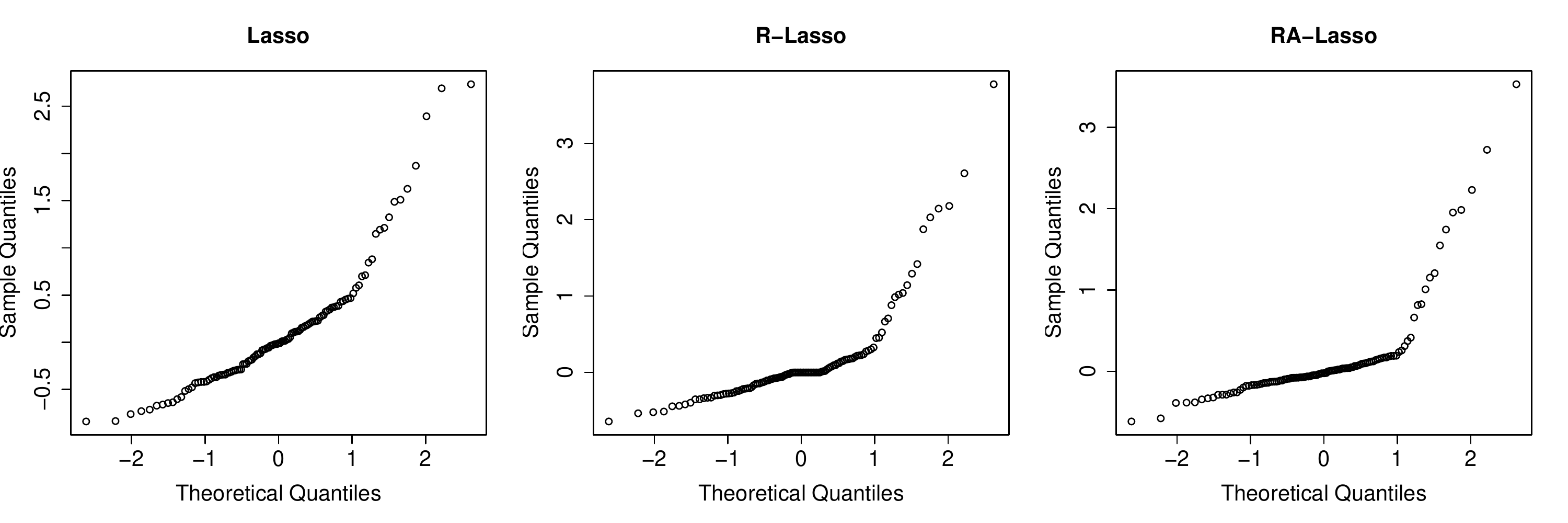}
  \caption{QQ plots of the residuals from three methods.\label{fig2}}
\end{figure}

In this section, we use a microarray data to illustrate the performance of Lasso, R-Lasso and
RA-Lasso. \cite{HUA11} studied the role of innate immune system on the development of
atherosclerosis by examining gene profiles from peripheral blood of 119 patients. The data were
collected using Illumina HumanRef8 V2.0 Bead Chip and are available on Gene Expression Omnibus. The
original study showed that the toll-like receptors (TLR) signaling pathway plays an important role
on triggering the innate immune system in face of atherosclerosis. Under this pathway, the ``TLR8''
gene was found to be a key atherosclerosis-associated gene. To further study the relationship
between this key gene and the other genes, we regressed it on another 464 genes from 12 different
pathways (TLR, CCC, CIR, IFNG, MAPK, RAPO, EXAPO, INAPO, DRS, NOD, EPO, CTR) that are related to
the TLR pathway. We applied Lasso, R-Lasso and RA-Lasso to this data. The tuning parameters for all
methods were chosen by using five-fold cross validation. Figure \ref{fig1} shows our choice of the
penalization parameter based on the cross validation results. For RA-Lasso, the choice of $\alpha$
was insensitive to the results and was fixed at 5. We then applied the three methods with the above
choice of tuning parameters to select significant genes. The QQ-plots of the residuals from the
three methods are shown in Figure \ref{fig2}. The selected genes by the three methods are reported
in Table \ref{tab4}. After the selection, we regressed the expression of TLR8 gene on the selected
genes, the $t$-values from the refittings are also reported in Table \ref{tab4}.


\begin{figure}[hbtp!]
  \centering
  \includegraphics[width=\textwidth]{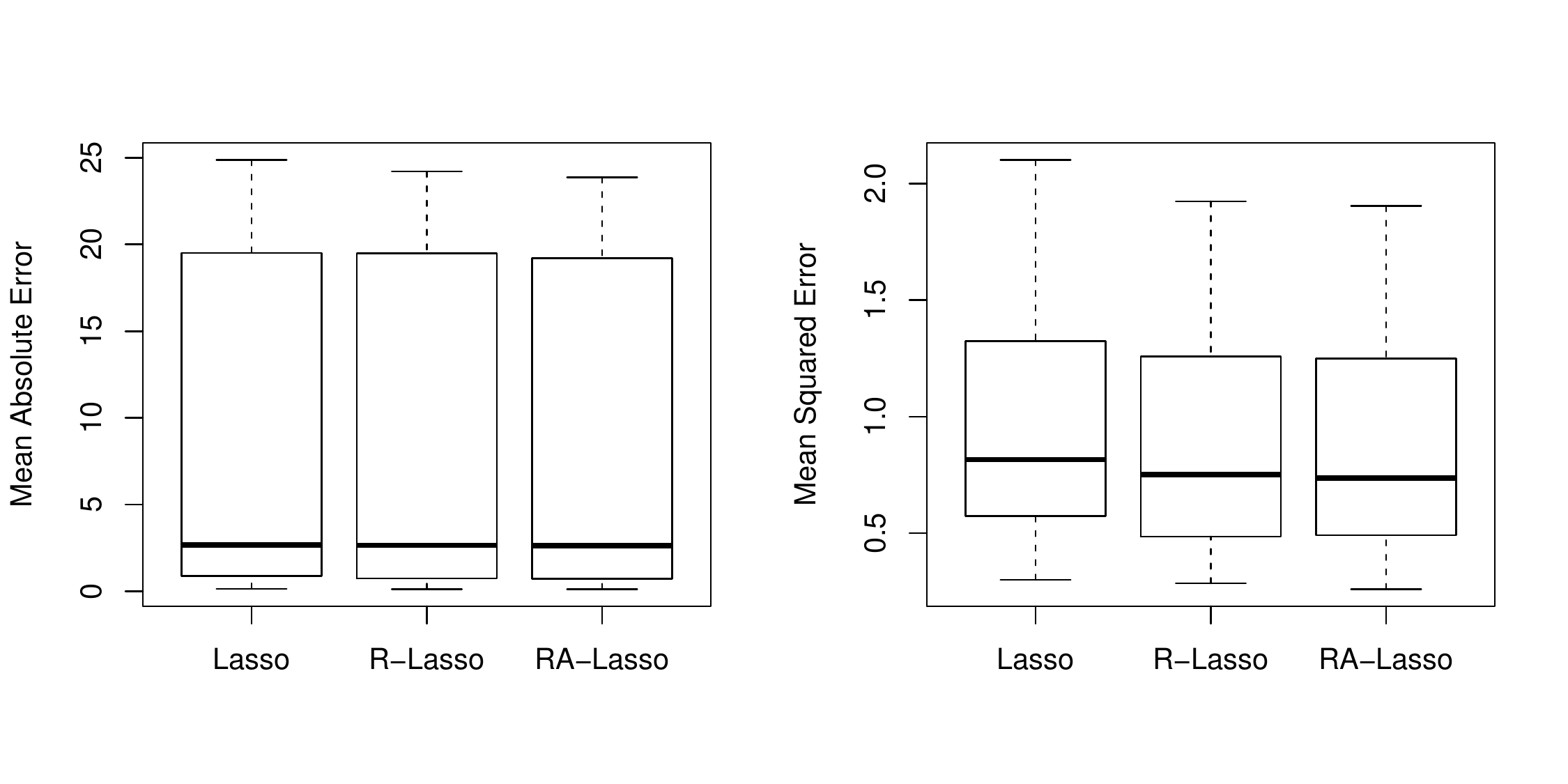}
  \vspace{-0.6in}
  \caption{Boxplot of the Mean Absolute/Squared Error of predictions.\label{fig3}}
\end{figure}

Table \ref{tab4} shows that Lasso only selected one gene. R-Lasso selected 17 genes. Our proposed
RA-Lasso selected 34 genes. Eight genes (CSF3, IL10, AKT1, TOLLIP, TLR1, SHC1, EPOR, and TJP1)
found by R-Lasso were also selected by RA-Lasso. Compared with Lasso and R-Lasso, our method
selected more genes, which could be useful for a second-stage confirmatory study.  It is clearly
seen from Figure \ref{fig2} that the residuals from the fitted regressions had heavy right tail and
skewed distribution. We know from the simulation studies in Section \ref{sec6} that RA-Lasso tends
to perform better than Lasso and R-Lasso in this situation. For further investigation, we randomly
chose 24 patients as the test set; applied three methods to the rest patients to obtain the
estimated coefficients, which in return were used to predict the responses of 24 patients. We
repeated the random splitting 100 times, the boxplots of the Mean Absolute/Squared Error of
predictions are shown in Figure \ref{fig3}. RA-Lasso has better predictions than Lasso and R-Lasso.

\begin{table}[hbtp!]
  \begin{tabular}[h]{l|cccccccc}
    \hline\hline
    Lasso    & CRK   &        &        &         &        &        & \\
    & 0.23  &        &        &         &        &        & \\
    \hline
    R-Lasso  & CSF3  & IL10   & AKT1   & KPNB1   & TLR2   & GRB2   & MAPK1 \\
    & -2.46 & 2.24   & 1.68   & 1.49    & 1.41   & -1.06  & 0.98 \\
    & DAPK2 & TOLLIP & TLR1   & TLR3    & SHC1   & PSMD1  & F12 \\
    & 0.7   & -0.68  & 0.52   & 0.33    & -0.28  & 0.27   & 0.24 \\
    & EPOR  & TJP1   & GAB2   &         &        &        & \\
    & -0.17 & -0.12  & -0.01  &         &        &        & \\
    \hline
    RA-Lasso & CSF3  & CD3E   & BTK    & CLSPN   & RELA   & AKT1   & IRS2 \\
    & -2.95 & 2.67   & 2.37   & 1.93    & 1.88   & 1.61   & 1.55 \\
    & IL10  & MAP2K4 & PMAIP1 & BCL2L11 & AKT3   & DUSP10 & IRF4 \\
    & 1.52  & 1.17   & -1.14  & -1.13   & -1.01  & 0.97   & -0.95 \\
    & IFI6  & TLR1   & PSMB8  & KPNB1   & IFNG   & FADD   & TJP1 \\
    & 0.86  & 0.82   & 0.79   & 0.77    & -0.74  & 0.65   & -0.57 \\
    & CR2   & IL2    & PSMC2  & HSPA8   & SHC1   & SPI1   & IFNA6 \\
    & 0.57  & -0.47  & 0.38   & -0.35   & -0.33  & -0.28  & 0.28 \\
    & FYN   & EPOR   & MASP1  & PRKCZ   & TOLLIP & BAK1   & \\
    & -0.24 & 0.24   & -0.24  & 0.24    & -0.19  & 0.14   & \\
    \hline\hline
  \end{tabular}
  \caption{Selected genes by Lasso, R-Lasso and RA-Lasso.  \label{tab4}}
\end{table}

\section{Proofs}
\label{sec:proofs}

\begin{proof}[\textbf{Proof of Theorem \ref{thm1}.}]
  Let $\ell(x)=x^2$. Since $\betas$ minimizes $\E\ell(y-\bx^T\bbeta)$, it follows from
  condition(C2) that
  \begin{equation}
    \label{eq7.1}
    \E[\ell(y-\bx^T\betaa)-\ell(y-\bx^T\betas)]=(\betaa-\betas)^T\E(\bx\bx^T)(\betaa-\betas)\geq
    \kappa_l  \ltwonorm{\betaa-\betas}^2.
  \end{equation}
  Let $g_{\alpha}(x)=\ell(x)-\ell_{\alpha}(x)=(|x|-\alpha^{-1})^2I(|x|>\alpha^{-1})$.  Then, since $\betas_\alpha$ is the minimizer of $E \ell_{\alpha}(y-\bx^T\bbeta)$, we have
  \begin{eqnarray*}
    &   & \E[\ell(y-\bx^T\betaa)-\ell(y-\bx^T\betas)]\\
    & = & \E[\ell(y-\bx^T\betaa)-\ell_{\alpha}(y-\bx^T\betaa)]
    +\E[\ell_{\alpha}(y-\bx^T\betaa)-\ell_{\alpha}(y-\bx^T\betas)]\\
    &  & \hspace{3ex}+\E[\ell_{\alpha}(y-\bx^T\betas)-\ell(y-\bx^T\betas)]\\
    & \leq & \E[g_{\alpha}(y-\bx^T\betaa)]-\E[g_{\alpha}(y-\bx^T\betas)].
  \end{eqnarray*}
  By Taylor's expansion, we have
  \begin{equation} \label{eq7.2}
    \E[\ell(y-\bx^T\betaa)-\ell_{\alpha}(y-\bx^T\betaa)]
    \leq  2\E[ (z-\alpha^{-1}) I(z>\alpha^{-1}) |\bx^T(\betaa-\betas)|],
  \end{equation}
  where $\tilde{\bbeta}$ is a vector lying between $\betas$ and $\betaa$ and $z=|y-\bx^T\tilde{\bbeta}|$.  With
  $\E_{\epsilon}$ denoting the conditional expectation with respect to $\epsilon$ given $\bx$, we have
  \begin{align*}
    \E_{\epsilon}[ (z-\alpha^{-1})I(z >\alpha^{-1})] \leq \E_{\epsilon}[ z I(z >\alpha^{-1})] \leq
    \alpha^{k-1}\E_{\epsilon} z^k.
  \end{align*}
  Therefore, $\E[\ell(y-\bx^T\betaa)-\ell(y-\bx^T\betas)]$ is further bounded by
  \begin{equation}
    \label{eq7.3}
    2\alpha^{k-1}\E(|y-\bx^T\tilde{\bbeta}|^k|\bx^T(\betaa-\betas)|)\leq
    2(2\alpha)^{k-1} \E[(M_k+|\bx^T(\tilde{\bbeta}-\betas)|^k)|\bx^T(\betaa-\betas)|],
  \end{equation}
  where the constant $M_k$ is defined in Condition (C1).
  Next, we show that
  $\lambda_{\max}(\E[(M_k+|\bx^T(\tilde{\bbeta}-\betas)|^k)^2\bx\bx^T])=O(1)$. Let $\bnu$ be a
  $p$-dimensional vector with $\ltwonorm{\bnu}=1$. By the Cauchy-Schwartz inequality,
  \begin{equation*}
    \E[(M_k+|\bx^T(\tilde{\bbeta}-\betas)|^k)^2(\bx^T\bnu)^2]\leq [\E(M_k+|\bx^T(\tilde{\bbeta}-\betas)|^k)^4]^{1/2}
    [\E(\bx^T\bnu)^4]^{1/2}.
  \end{equation*}
  By (C3), $\bx^T\bnu$ is sub-Gaussian with parameter $\kappa_0^2$. Under the assumption that
  $\betas$ and $\betaa$ are interior points of an $L_2$-ball with sufficiently large radius,
  $\bx^T(\tilde{\bbeta}-\betas)$ is sub-Gaussian with parameter $\kappa_0^2
  \ltwonorm{\tilde{\bbeta}-\betas}^2$, which is no larger than $\kappa_0^2
  \ltwonorm{\betaa-\betas}^2=O(1)$. Using the moment results of sub-Gaussian random variables
  \citep{RIV12}, $\E(\bx^T\bnu)^4 \leq 16\kappa_0^4=O(1)$. Similarly,
  $\E|\bx^T(\tilde{\bbeta}-\betas)|^{4k}\leq \E|\bx^T(\betas-\betaa)|^{4k}=O(1)$. Therefore,
  $\E[(M_k+|\bx^T(\tilde{\bbeta}-\betas)|^k)^4]=O(1)$. Hence, by definition,
  $\lambda_{\max}(\E[(M_k+|\bx^T(\tilde{\bbeta}-\betas)|^k)^2\bx\bx^T])=O(1)$. Using this result
  and (\ref{eq7.3}),
  \begin{align*}
    \E[\ell(y-\bx^T\betaa)-\ell(y-\bx^T\betas)]&\leq
    2(2\alpha)^{k-1}[\lambda_{\max}(\E[(M_k+|\bx^T(\tilde{\bbeta}-\betas)|^k)^2\bx\bx^T])]^{1/2}
    \ltwonorm{\betaa-\betas}\\
    &=O(\alpha^{k-1} \ltwonorm{\betaa-\betas}).
  \end{align*}
  This together with (\ref{eq7.1}) completes the proof.
\end{proof}

\begin{proof}[\textbf{Proof of Lemma \ref{lem1}.}]
  First of all, it follows from Lemma 1 of \cite{NRW12} that $\hat{\bDelta}=\hat{\bbeta}-\betaa\in
  \bbC_{\alpha \eta}$ on the event $\{\lambda_n\geq 2\supnorm{\nabla\cL_n(\betaa)} \}$.  Hence, we
  need to show that the event $\{\lambda_n\geq 2\supnorm{\nabla\cL_n(\betaa)} \}$ holds with high
  probability.  The latter will be established by using Bernstein's inequality along with the union
  bound.

  The gradient of $\cL_n$,
  \begin{equation}
    \label{eq7.4}
    \nabla \cL_n(\betaa)=\frac{1}{n} \sum_{i=1}^n \frac{2}{\alpha} \psi[\alpha(y_i-\bx_i^T\betaa)]\bx_i,
  \end{equation}
  where $\psi(x)=x$, for $|x|\leq 1$; $\psi(x)=1$, for $x>1$; and $\psi(x)=-1$, for $x<-1$.  Using $\alpha^{-1}|\psi(\alpha x)|\leq |x|$ and assumption (C3), we have
  \begin{align*}
    \E\{2\alpha^{-1}\psi[\alpha(y_i-\bx_i^T\betaa)]x_{ij} \}^2
    &\leq 4\E\{(y_i-\bx_i^T\betaa)^2x_{ij}^2 \}\\
    &\leq 8\E\{(\epsilon_i^2+|\bx_i^T(\betaa-\betas)|^2)x_{ij}^2\}\\
    &\leq v,
  \end{align*}
  where $v$ is a constant depending on $\kappa_0$ and $M_2$ and the last inequality follows from a
  similar argument as in the proof of Theorem \ref{thm1}. By (C3) and that $|\psi(x)|\leq 1$,
  $\psi[\alpha(y_i-\bx_i^T\betaa)]x_{ij}$ is also sub-Gaussian. For any $k\geq 3$, using the
  relation between the $k$th moment and the second moment of sub-Gaussian random variables
  \citep{RIV12},
  $$
  \E|\psi[\alpha(y_i-\bx_i^T\betaa)]x_{ij}|^k\leq \frac{k!}{2}L^{k-2}\E|\psi[\alpha(y_i-\bx_i^T\betaa)]x_{ij}|^2,
  $$
  where $L$ is a constant depending on $\kappa_0$ only.
  Hence,
  \begin{equation*}
    \E|2\alpha^{-1}\psi[\alpha(y_i-\bx_i^T\betaa)]x_{ij}|^k
    \leq \frac{k!}{2}(2L/\alpha)^{k-2} v.
  \end{equation*}
  By Bernstein inequality (Proposition 2.9 of \cite{MP07}) and note that
  $\E(\frac{2}{\alpha}\psi[\alpha(y_i-\bx_i^T\betaa)]\bx_i)=\bfsym{0}$, we have
  \begin{equation*}
    P\left(\left|\frac{1}{n} \sum_{i=1}^n
        \frac{2}{\alpha}\psi[\alpha(y_i-\bx_i^T\betaa)]x_{ij}\right| \geq
      \sqrt{\frac{2vt}{n}}+\frac{Lt}{\alpha n} \right)\leq 2\exp(-t).
  \end{equation*}
  Let $t=n\lambda_n^2/(32v)$ and observe that $\frac{2Lt}{\alpha n}\leq \sqrt{\frac{2vt}{n}}$ by the
  choice of $\lambda_n$ and $\alpha$. We have
  \begin{equation*}
    P\left(\left|\frac{1}{n} \sum_{i=1}^n
        \frac{2}{\alpha}\psi[\alpha(y_i-\bx_i^T\betaa)]x_{ij}\right| \geq
      \frac{\lambda_n}{2}\right) \leq 2\exp \left(-\frac{n\lambda_n^2}{32 v} \right).
  \end{equation*}
  It then follows from union inequality that
  \begin{align*}
    P\left (\supnorm{\frac{1}{n} \sum_{i=1}^n \frac{2}{\alpha}
        \psi[\alpha(y_i-\bx_i^T\betaa)]\bx_i}>\frac{\lambda_n}{2} \right )&\leq 2 \exp
    \left(-\frac{n\lambda_n^2}{32 v}+\log p \right)=2p^{-c_0},
  \end{align*}
  where $c_0=\kappa_{\lambda}^2/(32v)-1$. This completes the proof.
\end{proof}

\begin{proof}[\textbf{Proof of Lemma \ref{lem2}.}]
  Denote $\cL_n(\bbeta)=\frac{1}{n} \sum_{i=1}^n \la(y_i-\bx_i^T\bbeta)$. Applying a
  second-order Taylor expansion to $\cL_n(\bbeta)$ between $\betaa$ and $\betaa+\bDelta$, we conclude that for some $v\in [0,1]$,
  \begin{eqnarray}
    \delta\cL_n(\bDelta,\betaa)
    & = & \frac{1}{n} \sum_{i=1}^n
    \psi'[\alpha(y_i-\bx_i^T\betaa+v\bx_i^T\bDelta)](\bx_i^T\bDelta)^2, \label{eq7.5}
  \end{eqnarray}
  where $\psi'(x)=1$ for $|x|\leq 1$, and $\psi'(x)=0$ otherwise. Note that each term in
  \eqref{eq7.5} is nonnegative.  However, the quadratic component in \eqref{eq7.5} is not Lipschitz
  continuous with a bounded Lipschitz coefficient.  In order to apply the contraction theorem of
  \cite{LT91}, we introduce a truncation function that is Lipschitz and bound \eqref{eq7.5} from
  below.  Let
  \begin{equation}
    \label{eq7.6}
    \varphi_t(u)= u^2 I(|u|\leq t/2) +
    (t-u)^2 I(t/2\leq |u|\leq t),
  \end{equation}
  where $I(\cdot)$ is the indicator function.
  Clearly, $\varphi_t(u) \leq u^2$ and satisfies the Lipschitz condition with Lipschitz coefficient bounded by $2t$.  We first show
  \begin{equation}  \label{eq7.7}
    \delta\cL_n(\bDelta,\betaa)\geq \frac{1}{n} \sum_{i=1}^n \varphi_{\tau}(\bx_i^T\bDelta
    I(|y_i-\bx_i^T\betaa|\leq T)),
  \end{equation}
  for $0<\alpha \leq 1/(T+\tau)$, where the thresholds $T$ and $\tau$ will be chosen as in (\ref{eq7.10}).

  Let $A_i=\psi'[\alpha(y_i-\bx_i^T\betaa+v\bx_i^T\bDelta)](\bx_i^T\bDelta)^2$. We need only to
  show that
  $$
  A_i \geq \varphi_{\tau}(\bx_i^T\bDelta I(|y_i-\bx_i^T\betaa|\leq T)).
  $$
  When $|y_i-\bx_i^T\betaa|>T$ or $|\bx_i^T\bDelta|>\tau$, the right hand side is zero and the
  inequality holds trivially.  Thus, we need only to consider the case $|y_i-\bx_i^T\betaa|\leq T$
  and $|\bx_i^T\bDelta| \leq \tau$.  In this case,
  $$
  |\alpha(y_i-\bx_i^T\betaa+v\bx_i^T\bDelta)| \leq \alpha (T + \tau) \leq 1,
  $$
  and hence $\psi'[\alpha(y_i-\bx_i^T\betaa+v\bx_i^T\bDelta)]=1$.  Using this,
  $$
  A_i = (\bx_i^T\bDelta)^2 \geq \varphi_{\tau}(\bx_i^T\bDelta) = \varphi_{\tau}(\bx_i^T\bDelta
  I(|y_i-\bx_i^T\betaa|\leq T)).
  $$

  Using (\ref{eq7.7}), to prove the Lemma, we need to show that, for any fixed $\delta\in (0,1]$,
  with high probability
  \begin{equation}
    \label{eq7.8}
    \bbP_n\varphi_{\tau}(\bx^T\bDelta I(|y-\bx^T\betaa|\leq T))\geq \kappa_1 \ltwonorm{\bDelta}
    \{\ltwonorm{\bDelta}-\kappa_2\sqLogpN \lonenorm{\bDelta} \}, ~
    \text{for all } \ltwonorm{\bDelta}=\delta,
  \end{equation}
  where constants $\kappa_1$ and $\kappa_2$ do not depend on $\alpha$
  and
  $$
  \bbP_n\varphi_{\tau}(\bx^T\bDelta I(|y-\bx^T\betaa|\leq T))=\frac{1}{n} \sum_{i=1}^n
  \varphi_{\tau}(\bx_i^T\bDelta I(|y_i-\bx_i^T\betaa|\leq T)).
  $$
  This is equivalent to proving (\ref{eq7.8}) for $\delta=1$.  Indeed,
  from the definition (\ref{eq7.6}), for any $d>0$ and $z\in \mathbb{R}$, we have
  $\varphi_d(dz)=d^2\varphi_1(z)$. Thus, the event
  \begin{equation}\label{eq7.9}
    \bbP_n\varphi_{\tau_1 }(\bx^T\bDelta I(|y-\bx^T\betaa|\leq T))\geq \kappa_1
    \{1 -\kappa_2\sqLogpN \lonenorm{\bDelta} \}, \text{ for all } \ltwonorm{\bDelta}=1
  \end{equation}
  is the same as the event
  $$
  \bbP_n\varphi_{\tau_1}(\bx^T(\bDelta/\ltwonorm{\bDelta}) I(|y-\bx^T\betaa|\leq T))\geq \kappa_1
  \{1-\kappa_2\sqLogpN \lonenorm{\bDelta}/\ltwonorm{\bDelta} \},
  $$
  which equals to the event \eqref{eq7.8} with $\tau = \delta\tau_1$.

  To establish \eqref{eq7.9}, let us consider its complementary event.  Define
  \begin{equation*}
    \fd(\bx)=\bx^T\bDelta I(|y-\bx^T\betaa|\leq T) \qquad \text{and} \qquad
    \gd(\bx)=\varphi_{\tau}(\fd(\bx)).
  \end{equation*}
  Let $\bbS_2(1)$ be the unit sphere with $L_2$-radius one, and $\bbS_1(t)$ be the sphere of $L_1$-radius $t$, which is to be chosen later. The complementary event of \eqref{eq7.9} is given by
  \begin{equation*}
    \Big\{\bbP_n[\gd(\bx)]<\kappa_1\{1-\kappa_2 \sqLogpN \lonenorm{\bDelta}\}, \text{ for some }
    \bDelta\in \bbS_2(1)\Big\}.
  \end{equation*}
  Our goal is to show that the probability of this event is very small,  which is demonstrated
  through the following three steps.
  \begin{itemize}
  \item [(a)] First, we show that with the following choice of truncation
    \begin{equation}
      \label{eq7.10}
      T^2={1024\ko^4\kl^{-2}{M_k}^{{2}/{k}}} \qquad \text{ and } \qquad
      \tau^2=\max\{32\ko^2\log(12\kl^{-1}\ko^2), 1\},
    \end{equation}
    for any fixed $\bDelta\in \bbS_2(1)$, we have
    \begin{equation}
      \label{eq7.11}
      \E[\gd(\bx)]\geq \kl/2.
    \end{equation}

  \item [(b)] Second, with $Z(t)=\sidx |\bbP_n[\gd(\bx)]-\E[\gd(\bx)]|$,
    we prove the tail probability bound for $Z(t)$ is bounded by
    \begin{equation}
      \label{eq7.12}
      P(Z(t)\geq {\kl}/{4}+40\tau^2\ko t\sqLogpN)\leq \exp(-c_1'n-c_2't^2\log p),
    \end{equation}
    for each given $t$.

  \item [(c)] Finally, we use a standard peeling argument \citep{ALE87,VDG00} to establish
    $$
    P\Bigl \{\exists \bDelta\in \bbS_2(1):  Z(\lonenorm{\bDelta})\geq
    {\kl}/{4}+40\tau^2\ko \lonenorm{\bDelta}\sqLogpN \Bigr \} \leq  \exp(-c_1 n-c_2 \log p).
    $$
  \end{itemize}
  The result (c) together with (\ref{eq7.11}) show that the probability of the complementary event
  of (\ref{eq7.9}) with $\kappa_1=\kappa_{l}/4$ and $\kappa_2=40\tau^2\kappa_0\kappa_1^{-1}$ is
  bounded by $\exp(-c_1n-c_2\log p)$, which completes the proof.

  We first prove (\ref{eq7.11}). In fact, by condition (C2), for any $\bDelta\in
  \bbS_2(1)$, $\E[(\bx^T\bDelta)^2]\geq \kl \ltwonorm{\bDelta}^2=\kl$. So, it suffices to show that
  $\E[(\bx^T\bDelta)^2-\gd(\bx)]\leq \kl/2$.

  Note that, $\gd(\bx)=(\bx^T\bDelta)^2$ for all $\bx$ such that $|y-\bx^T\betaa|\leq T$ and
  $|\bx^T\bDelta|\leq \tau/2$. Therefore, we have
  \begin{equation}
    \label{eq7.13}
    \E[(\bx^T\bDelta)^2-\gd(\bx)]\leq \E[(\bx^T\bDelta)^2I(|y-\bx^T\betaa|\geq
    T)]+\E[(\bx^T\bDelta)^2I(|\bx^T\bDelta|\geq \tau/2)].
  \end{equation}
  To bound the first term on the right hand side of (\ref{eq7.13}), it follows from the
  Cauchy-Schwartz inequality that
  \begin{equation*}
    \E[(\bx^T\bDelta)^2I(|y-\bx^T\betaa|\geq T)] \leq [\E(\bx^T\bDelta)^4]^{1/2}
    [P(|y-\bx^T\betaa|\geq T)]^{1/2}.
  \end{equation*}
  Since $\bx^T\bDelta$ is sub-Gaussian with parameter at most $\ko^2$ by assumption (C3), we have
  $\E(\bx^T\bDelta)^4\leq 16\ko^4$. Meanwhile, for any $0<\alpha\leq 1/(T+\tau)$, it follows from the Chebyshev
  inequality and Theorem \ref{thm1} that
  \begin{align*}
    T^2 P(|y-\bx^T\betaa|\geq T) &\leq \E[(y-\bx^T\betaa)^2]\\
    &\leq 2\E \epsilon^2+2\E[\bx^T(\betas-\betaa)]^2\\
    &\leq 2 {M_k}^{2/k}+O(\alpha^{2k-2})\\
    &\leq 4 {M_k}^{2/k}.
  \end{align*}
  To bound the second term on the right hand side of (\ref{eq7.13}), by the concentration inequality
  of sub-Gaussian variables, we have
  \begin{equation*}
    P(|\bx^T\bDelta|\geq \tau/2)\leq 2\exp\{-{\tau^2}/({8\ko^2})\}.
  \end{equation*}
  Then, by the choice of $T$ and $\tau$ in (\ref{eq7.10}),
  \begin{equation*}
    \E[(\bx^T\bDelta)^2I(|y-\bx^T\betaa|\geq T)]\leq \frac{\kl}{4} \qquad \text{ and } \qquad
    \E[(\bx^T\bDelta)^2I(|\bx^T\bDelta|\geq \tau/2)] \leq \frac{\kl}{4}.
  \end{equation*}
  Hence, (\ref{eq7.11}) follows.

  Next, we give the tail bound as in (b). Indeed, for any $\bDelta\in \bbS_2(1)$, we have
  $\supnorm{\gd}\leq \tau^2$. Therefore, by Massart concentration inequality (Theorem 14.2 of
  \cite{BV11}), for any $z>0$, we have $P(Z(t)\geq \E Z(t)+z)\leq \exp(-\frac{nz^2}{32\tau^4})$. By
  choosing $z=\kl/4+16\tau^2\ko t\sqLogpN$, we have
  \begin{equation}
    \label{eq7.14}
    P(Z(t)\geq \E Z(t)+z)\leq  \exp \Bigl (-\frac{n\kl^2}{512\tau^4}-8\ko^2t^2\log p \Bigr ).
  \end{equation}

  Next, we bound $\E Z(t)$. Let $\{\omega_i \}_{i=1}^n$ be an i.i.d. sequence of Rademacher
  variables. A symmetrization theorem (Theorem 14.3 of \cite{BV11}) yields
  \begin{equation*}
    \E[Z(t)]\leq 2\E\left[\sidx |\frac{1}{n} \sum_{i=1}^n \omega_i
      \gd(\bx_i)|\right]=2\E\left[\sidx |\frac{1}{n}
      \sum_{i=1}^n \omega_i\varphi_{\tau}(\fd(\bx_i))|\right].
  \end{equation*}
  By definition, the function $\varphi_{\tau}$ is Lipschitz with parameter at most $2\tau
  \leq 2\tau^2$ and $\varphi_{\tau}(0)=0$. Therefore, by the Ledoux-Talagrand contraction theorem
  (\cite{LT91}, p.112), we have
  \begin{align*}
    \E[Z(t)]&\leq 8\tau^2\E\left[\sidx |\frac{1}{n} \sum_{i=1}^n \omega_i \fd(\bx_i)|\right] \\
    &= 8\tau^2\E\left[\sidx |\frac{1}{n} \sum_{i=1}^n \omega_i\bx_i^T\bDelta
      I(|y_i-\bx_i^T\betaa|\leq T)|\right]\\
    &\leq 8\tau^2t \E \supnorm{\frac{1}{n} \sum_{i=1}^n \omega_i\bx_iI(|y_i-\bx_i^T\betaa|\leq T)}.
  \end{align*}
  Since the variables $\{x_{ij} \}_{i=1}^n$ are zero-mean i.i.d. sub-Gaussian with parameter at most
  $\ko^2$, so are $\{\omega_ix_{ij}I(|y_i-\bx_i^T\betaa|\leq T) \}_{i=1}^n$. Since $\E
  \supnorm{\frac{1}{n} \sum_{i=1}^n \omega_i\bx_iI(|y_i-\bx_i^T\betaa|\leq T)}$ is the maxima of $p$
  such terms, known bounds on the expectation of sub-Gaussian maxima (e.g. see
  \cite{LT91}, p.79) yield
  \begin{equation*}
    \E \supnorm{\frac{1}{n} \sum_{i=1}^n \omega_i\bx_iI(|y_i-\bx_i^T\betaa|\leq T)} \leq 3\ko \sqLogpN.
  \end{equation*}
  Hence,
  \begin{equation}
    \label{eq7.15}
    \E[Z(t)]\leq 24\tau^2\ko t \sqLogpN.
  \end{equation}
  Combining (\ref{eq7.14}) and (\ref{eq7.15}), we have
  \begin{equation*}
    P\Big(Z(t)\geq {\kl}/{4}+40\tau^2\ko t\sqLogpN\Big)\leq \exp(-c_1'n-c_2't^2\log p),
  \end{equation*}
  where constants $c_1'$ and $c_2'$ depends on $\kl$ and $\ko$ only.  This result holds for each given $t$.

  Next, we furnish the peeling argument in (c). Let $h(\lonenorm{\bDelta})={\kl}/{8}+20\tau^2\ko
  \lonenorm{\bDelta}\sqLogpN$ and $B=\{\exists \bDelta\in \bbS_2(1): Z(\lonenorm{\bDelta})\geq
  2h(\lonenorm{\bDelta})\}$. Since $h(\lonenorm{\bDelta}) \geq \kappa_l/8$, the set can be covered
  by partition $\{B_m\}_{m=1}^\infty$ with $B_m=\{\bDelta\in \bbS_2(1): 2^{m-4}\kl\leq
  h(\lonenorm{\bDelta})\leq 2^{m-3} \kl \}$. Thus, by union bound,
  \begin{eqnarray*}
    P(B) & \leq & \sum_{m=1}^{\infty} P(\bDelta\in B_m \text{ such that }  Z(\lonenorm{\bDelta})\geq
    2h(\lonenorm{\bDelta}))\\
    & \leq & \sum_{m=1}^{\infty} P( Z(\lonenorm{\bDelta})\geq
    2^{m-3}\kl)
  \end{eqnarray*}
  since $Z(\lonenorm{\bDelta})\geq 2^{m-3}\kl$ for $\bDelta\in B_m$.
  By letting $2^{m-3}\kl={\kl}/{4}+40\tau^2\ko t\sqLogpN$ as in (\ref{eq7.12}) and solving for $t$, by \eqref{eq7.12}, we obtain
  \begin{align*}
    P(B) %
    &\leq
    \sum_{m=1}^{\infty}\exp\left(-c_1'n-\frac{c_2'\kl^2(2^{m-1}-1)^2n}{\tau^4\ko^2}\right)\\
    &\leq \exp(-c_1'n)+
    \sum_{m=2}^{\infty}\exp\left(-c_1'n-\frac{c_2'n\kl^22^{2m-4}}{\tau^4\ko^2}\right)\\
    &\leq c_1\exp(-c_2n),
  \end{align*}
  where the last inequality follows from sum of geometric series.
\end{proof}

\begin{proof}[\textbf{Proof of Lemma \ref{lem3}.}]
  Note that,
  \begin{equation}
    \label{eq7.16}
    \Rq\geq \sum_{j=1}^p |\beta_{\alpha,j}^{\ast}|^q \geq \sum_{j\in S_{\alpha \eta}}
    |\beta_{\alpha,j}^{\ast}|^q\geq \eta^q|S_{\alpha \eta}|.
  \end{equation}
  Therefore, $|S_{\alpha \eta}|\leq \eta^{-q}\Rq$. Let $S_{\alpha
    \eta}^c=\{1,2,\ldots,p \}\backslash S_{\alpha \eta}$, we have
  \begin{equation}
    \label{eq7.17}
    \lonenorm{\bbeta_{S_{\alpha \eta}^c}^{\ast}}=\sum_{j\in S_{\alpha\eta}^c}|\beta_{\alpha,j}^{\ast}| =
    \sum_{j\in S_{\alpha \eta}^c} |\beta_{\alpha,j}^{\ast}|^q|\beta_{\alpha,j}^{\ast}|^{1-q}\leq
    \Rq\eta^{1-q}.
  \end{equation}
  Hence, for any $\bDelta\in \bbC_{\alpha \eta}$, we have
  \begin{align*}
    \lonenorm{\bDelta}= \lonenorm{\bDelta_{S_{\alpha \eta}}}+  \lonenorm{\bDelta_{S_{\alpha
          \eta}^c}} \leq 4 \lonenorm{\bDelta_{S_{\alpha \eta}}}+4 \lonenorm{\bbeta^{\ast}_{S_{\alpha
          \eta}^c}}.
  \end{align*}
  By the Cauchy-Schwartz inequality and \eqref{eq7.17}, we can bound further that
  \begin{align*}
    \lonenorm{\bDelta} &\leq 4 \sqrt{|S_{\alpha \eta}|} \ltwonorm{\bDelta}+4\Rq\eta^{1-q} \leq
    4\Rq^{1/2}\eta^{-q/2}\ltwonorm{\bDelta}+4\Rq\eta^{1-q}.
  \end{align*}
  It then follows from Lemma \ref{lem2} that
  \begin{align*}
    \delta\cL_n(\bDelta,\betaa)&\geq \kappa_1
    \ltwonorm{\bDelta}\{\ltwonorm{\bDelta}-\kappa_2\sqLogpN
    [4\Rq^{1/2}\eta^{-q/2}\ltwonorm{\bDelta}+4\Rq\eta^{1-q}] \}\\
    &=\Big(\kappa_1-4\kappa_1\kappa_2\Rq^{1/2}\eta^{-q/2}\sqLogpN\Big)
    \ltwonorm{\bDelta}^2-4\kappa_2\Rq\eta^{1-q}\sqLogpN.
  \end{align*}
  With $\lambda_n=\kappa_{\lambda} \sqLogpN$ and $\eta=\lambda_n$, it
  holds that
  \begin{equation*}
    4\kappa_1\kappa_2\Rq^{1/2}\eta^{-q/2}\sqrt{\frac{\log p}{n}}
    =4\kappa_1\kappa_2\Rq^{1/2}{\kappa_{\lambda}}^{-q/2}\left(\frac{\log p}{n} \right)^{(1-q)/2},
  \end{equation*}
  which is no larger than $\kappa_1/2$ under assumption (\ref{eq2.7}). On
  the other hand,
  \begin{equation*}
    4\Rq\kappa_2\eta^{1-q}\sqrt{\frac{\log p}{n}}
    =4\Rq\kappa_2{\kappa_{\lambda}}^{1-q} \left(\frac{\log p}{n} \right)^{1-(q/2)}.
  \end{equation*}
  Therefore, RSC holds with $\kappa_{\cL}=\frac{\kappa_1}{2}$ and $\tau^2_{\cL}=4\Rq\kappa_2
  {\kappa_{\lambda}}^{1-q}(\frac{\log p}{n})^{1-(q/2)}$.
\end{proof}

\begin{proof}[\textbf{Proof of Theorem \ref{thm2}.}]
  Let $A_1$ and $A_2$ denote the events that Lemma \ref{lem1} and Lemma \ref{lem3} hold,
  respectively. By Theorem 1 of \cite{NRW12}, within $A_1\cap A_2$, it holds that
  \begin{align*}
    \ltwonorm{\bDelta}^2&\leq 9 \frac{\lambda_n^2}{\kappa_{\cL}^2}|S_{\alpha
      \eta}|+\frac{\lambda_n}{\kappa_{\cL}^2}\{2\tau_{\cL}^2
    + 4\lonenorm{\bbeta^{\ast}_{S_{\alpha \eta}^c}} \}\\
    &\leq \frac{36\lambda_n^2\Rq}{\kappa_1^2\eta^q}+\frac{4\lambda_n}{\kappa_1^2} \left\{8\Rq\kappa_2
      {\kappa_{\lambda}}^{1-q}\left(\frac{\log p}{n}
      \right)^{1-(q/2)}+4\Rq\eta^{1-q} \right\}\\
    &\overset{(i)}{=}\frac{36}{\kappa_1^2}\Rq\lambda_n^{2-q}+\frac{16}{\kappa_1^2} \Rq\lambda_n^{2-q}
    \left\{2\kappa_2\left(\frac{\log p}{n}\right)^{\frac{1-q}{2}}+1\right\}  \\
    &=O\left(R_q\lambda_n^{2-q}\right) =O\left(R_q[(\log p)/n]^{1-(q/2)}\right),
  \end{align*}
  where (i) follows from the choice of $\eta=\lambda_n$. On the other hand, by Lemma \ref{lem1}
  and \ref{lem3}, $P(A_1\cap A_2)\geq 1-2p^{-c_0}-c_1\exp(-c_2n)$.
\end{proof}

\begin{proof}[\textbf{Proof of Lemma \ref{lem4}.}]
  From the proof of Lemma \ref{lem2}, we can see that \eqref{eq2.6} indeed holds for all $\bbeta$
  and $\bDelta\in \{\bDelta:\ltwonorm{\bDelta}\leq 1 \}$ that
  \begin{align*}
    \delta\cL_n(\bDelta,\bbeta)&\geq \kappa_1 \ltwonorm{\bDelta}^2-\kappa_1\kappa_2 \ltwonorm{\bDelta}
    \lonenorm{\bDelta} \sqLogpN.
  \end{align*}
  Using the fact that $ab \leq (a^2+b^2)/2$, we conclude that
  \begin{align*}
    \delta\cL_n(\bDelta,\bbeta)&\geq \kappa_1 \ltwonorm{\bDelta}^2-\left(\frac{1}{2}\kappa_1
      \ltwonorm{\bDelta}^2+\frac{1}{2} \kappa_1\kappa_2^2 \lonenorm{\bDelta}^2 \left(\frac{\log
          p}{n}\right)\right).
  \end{align*}
  Therefore, (\ref{eq3.2}) holds with $\gamma_l=\kappa_1$ and $\tau_l=\kappa_1\kappa_2^2(\log p)/(2n)$. Meanwhile, since
  $|\psi'(\cdot)|\leq 1$, it follows from (\ref{eq7.5}) that
  \begin{equation*}
    \delta\cL_n(\bDelta,\bbeta)\leq \frac{1}{n} \sum_{i=1}^n(\bx_i^T\bDelta)^2.
  \end{equation*}
  Under the sub-Gaussianity assumption (C3), it follows from some existing
  work (e.g. page 18 of \cite{LWAN13}) that, with probability great than
  $1-c_1\exp(-c_2n)$, it holds that
  \begin{equation*}
    \frac{1}{n} \sum_{i=1}^n (\bx_i^T\bDelta)^2\leq \ku\left(\frac{3}{2}\ltwonorm{\bDelta}^2
      +\frac{\log p}{n} \lonenorm{\bDelta}^2\right),
  \end{equation*}
  where $c_1$ and $c_2$ are some generic constants. Hence, (\ref{eq3.3}) holds with
  $\gamma_u={3\ku}$ and $\tau_u=\ku (\log p)/n$.
\end{proof}

\begin{proof}[\textbf{Proof of Theorem \ref{thm4}.}]
  We prove the theorem by the following two steps:

  (a) We first show that, for any $\delta^2\geq \varepsilon^2/(1-\kappa)$,
  $\phi(\hat{\bbeta}^t)-\phi(\hat{\bbeta})\leq \delta^2$, for all $t$ greater than the right hand
  side of (\ref{eq7.20}), where $\kappa\in [0,1)$ is a contraction constant and $\varepsilon$ is a
  tolerance parameter, which will be given in (\ref{eq7.21}) and (\ref{eq7.22}), respectively.

  (b) We use RSC condition (\ref{eq3.2}) to transform the upper bound of
  $\phi(\hat{\bbeta}^t)-\phi(\hat{\bbeta})$ into the upper bound of
  $\ltwonorm{\hat{\bbeta}^t-\hat{\bbeta}} $.

  For step (a), since our loss function is convex, we apply Theorem 2 of \cite{ANW12}. In
  order for our proof to be self-contained, we cite their theorem as the follows:

  \emph{[Theorem 2 of \cite{ANW12}] Suppose for any data set $Z_1^n$, the loss function
    $\cL_n(.;Z_1^n)$ is convex and differentiable and the regularizer $\cR$ is a norm. Consider the
    optimization problem of $\hat{\theta}=\argmin_{\cR(\theta)\leq \rho}\{\cL_n(\theta;
    Z_1^n)+\lambda_n \cR(\theta) \}$ for a radius $\rho$ such that $\theta^{\ast}$ is feasible,
    where $\theta^{\ast}=\argmin \E\cL_n(\theta;Z_1^n)$, and a regularization parameter $\lambda_n$
    satisfying bound
    \begin{equation}
      \label{eq7.18}
      \lambda_n\geq 2\cR^{\ast}(\nabla \cL_n(\theta^{\ast})),
    \end{equation}
    where $\cR^{\ast}$ is the dual norm of the regularizer. In addition, suppose that the loss
    function $\cL_n$ satisfies the RSC/RSM condition with parameters ($\gamma_l$, $\tau_l$) and
    ($\gamma_u$, $\tau_u$), respectively. Let $(\mathcal{M},\bar{\mathcal{M}}^{\perp})$ be any
    $\cR$-decomposable pair of subspaces such that
    \begin{equation}
      \label{eq7.19}
      \kappa = \left\{1-\frac{\bar{\gamma}_l}{4\gamma_u}
        +\frac{64\Psi^2(\bar{\mathcal{M}})\tau_u}{\bar{\gamma}_l} \right\}\xi \in [0,1)
      \qquad \text{ and } \qquad
      \frac{32 \rho}{1-\kappa}\xi\chi \leq \lambda_n,
    \end{equation}
    where $\Psi(\bar{\mathcal{M}})=\sup_{\theta\in \bar{\mathcal{M}}\backslash \{0\}}
    {\cR(\theta)}/{\ltwonorm{\theta}}$,
    $\bar{\gamma}_l=\gamma_l-64\tau_l\Psi^2(\bar{\mathcal{M}})$, $\xi=
    (1-{64\tau_u\bar{\gamma}_l^{-1}\Psi^2(\bar{\mathcal{M}})})^{-1}$, and
    $\chi=2\left({\bar{\gamma}_l}/(4\gamma_u)
      +128\tau_u\bar{\gamma}_l^{-1}\Psi^2(\bar{\mathcal{M}})\right)\tau_l+8\tau_u+2\tau_l$. Denote
    $\varepsilon^2= 8\xi\chi\left(6\Psi(\bar{\mathcal{M}})\ltwonorm{\hat{\theta}-\theta^{\ast}}
      +8\cR(\Pi_{\mathcal{M}^{\perp}}(\theta^{\ast}))\right)^2$, where
    $\Pi_{\mathcal{M}^{\perp}}(\theta^{\ast})$ is the projection of $\theta^{\ast}$ onto
    $\mathcal{M}^{\perp}$. Then for any $\delta^2\geq \varepsilon^2/(1-\kappa)$, we have
    $\phi_n(\hat{\theta}^t)-\phi_n(\hat{\theta})\leq \delta^2$ for all
    \begin{equation}
      \label{eq7.20}
      t\geq \frac{2\log((\phi_n(\theta^0)-\phi_n(\hat{\theta}))/\delta^2)}{\log(1/\kappa)}+
      \log_2\log_2 \left(\frac{\rho\lambda_n}{\delta^2} \right)\left(1+\frac{\log2}{\log(1/\kappa)}
      \right),
    \end{equation}
    where $\phi_n(\theta)= \cL_n(\theta;Z_1^n)+\lambda_n\cR(\theta)$, $\hat{\theta}^t$ is the
    solution by the gradient descent algorithm after $t^{\text{th}}$ iteration, and $\theta^0$ is
    the initial value of $\theta$.
  }

  In fact, Theorem 2 of \cite{ANW12} is a deterministic statement for all choices of pairs
  $(\mathcal{M},\bar{\mathcal{M}}^{\perp})$. From Lemma \ref{lem1} and Lemma \ref{lem4}, we have
  shown that with our choice of $\lambda_n$, the RA-quadratic loss function satisfy (\ref{eq7.18})
  and RSC/RSM with probability at least $1-2p^{-c_0}$ and $1-c_1\exp(-c_2n)$, respectively. Hence,
  Theorem 2 of \cite{ANW12} applies to our problem with high probability. We further choose the
  pair $(\mathcal{M},\bar{\mathcal{M}}^{\perp})=(S_{\alpha\eta},S_{\alpha\eta}^c)$ and give the
  explicit expression of constants for our problem as the follows:
  \begin{align}
    \kappa&=\left\{1-\frac{\bar{\gamma}_l}{4\gamma_u}+\frac{64\ku|S_{\alpha \eta}| \frac{\log
          p}{n}}{\bar{\gamma}_l} \right\}\left(1-\frac{64\kappa_u|S_{\alpha \eta}|\frac{\log
          p}{n}}{\bar{\gamma}_l}\right)^{-1} \label{eq7.21}\\
    \varepsilon^2&=8\xi\chi\left(6 \sqrt{|S_{\alpha \eta}|} \ltwonorm{\hat{\bbeta}-\betaa}+8
      \lonenorm{\betas_{S_{\alpha \eta}^c}}\right)^2, \label{eq7.22}
  \end{align}
  where $\bar{\gamma}_l=\kappa_1-32\kappa_1\kappa_2^2|S_{\alpha \eta}|{(\log p)}/{n}$,
  $\xi=\{1-64\ku |S_{\alpha \eta}|(\log p)/(n\bar{\gamma}_l)\}^{-1}$, and
  $\chi=2\{{\bar{\gamma}_l}/(4\gamma_u)+{128\tau_u|S_{\alpha
      \eta}|}/{\bar{\gamma}_l}+1 \}\tau_l+8\tau_u$. It remains to check (\ref{eq7.19}). By
  (\ref{eq7.21}), $\kappa\in [0,1)$ is equivalent to requiring
  \begin{equation}
    \label{eq7.23}
    |S_{\alpha\eta}|\frac{\log p}{n}< \frac{\bar{\gamma}_l^2}{1536\ku^2}
  \end{equation}
  With $\eta=\lambda_n$, it follows from (\ref{eq7.16}) that
  \begin{equation*}
    |S_{\alpha\eta}|\frac{\log p}{n} \leq \Rq\eta^{-q}\frac{\log p}{n}\leq
    \kappa_{\lambda}^{-q}\Rq\left(\frac{\log p}{n}\right)^{1-(q/2)}.
  \end{equation*}
  Hence, (\ref{eq7.23}) holds when $n$ is sufficiently large. Moreover, from \eqref{eq7.19} we
  need
  \begin{equation*}
    \lambda_n\geq \frac{32\rho}{1-\kappa}\left(1-\frac{64\ku|S_{\alpha\eta}|\frac{\log
          p}{n}}{\bar{\gamma}_l}\right)^{-1} \left[1+\kappa_1\kappa_2^2
      \left(\frac{\bar{\gamma}_l}{12\ku}+\frac{128\ku|S_{\alpha\eta}|\frac{\log
            p}{n}}{\bar{\gamma}_l} \right)+8\ku\right]\frac{\log p}{n},
  \end{equation*}
  which is satisfied under the stated assumption. It then follows from Theorem 2 of
  \cite{ANW12} that, for any $\delta^2\geq {\varepsilon^2}/(1-\kappa)$,
  $\phi(\hat{\bbeta}^t)-\phi(\hat{\bbeta})\leq \delta^2$, for all iterations $t$ greater than the
  right hand side of (\ref{eq7.20}).

  For step (b), it follows from RSC condition that
  \begin{equation*}
    \cL_n(\hat{\bbeta}^t)-\cL_n(\hat{\bbeta})-[\nabla\cL_n(\hat{\bbeta})]^T
    (\hat{\bbeta}^t-\hat{\bbeta}) \geq \frac{\gamma_l}{2}
    \ltwonorm{\hat{\bbeta}^t-\hat{\bbeta}}^2-\tau_l
    \lonenorm{\hat{\bbeta}^t-\hat{\bbeta}}^2.
  \end{equation*}
  Then we have
  \begin{align*}
    \phi(\hat{\bbeta}^t)-\phi(\hat{\bbeta})&=\cL_n(\hat{\bbeta}^t)-\cL_n(\hat{\bbeta})+
    \lambda_n(\lonenorm{\hat{\bbeta}^t}-\lonenorm{\hat{\bbeta}})\\
    &\geq [\nabla\cL_n(\hat{\bbeta})]^T(\hat{\bbeta}^t-\hat{\bbeta})+
    \lambda_n(\lonenorm{\hat{\bbeta}^t}-\lonenorm{\hat{\bbeta}})+\frac{\gamma_l}{2}
    \ltwonorm{\hat{\bbeta}^t-\hat{\bbeta}}^2-\tau_l \lonenorm{\hat{\bbeta}^t-\hat{\bbeta}}^2.
  \end{align*}
  Since $\hat{\bbeta}$ is the minimizer of $\phi(\bbeta)$, by the first-order condition, $[\nabla
  \cL_n(\hat{\bbeta})+\lambda_n\nabla \lonenorm{\hat{\bbeta}}]^T(\hat{\bbeta}^t-\hat{\bbeta})\geq
  0$. Therefore,
  \begin{equation*}
    \phi(\hat{\bbeta}^t)-\phi(\hat{\bbeta})\geq -\lambda_n[\nabla
    \lonenorm{\hat{\bbeta}}]^T(\hat{\bbeta}^t-\hat{\bbeta})+
    \lambda_n(\lonenorm{\hat{\bbeta}^t}-\lonenorm{\hat{\bbeta}})+\frac{\gamma_l}{2}
    \ltwonorm{\hat{\bbeta}^t-\hat{\bbeta}}^2-\tau_l \lonenorm{\hat{\bbeta}^t-\hat{\bbeta}}^2.
  \end{equation*}
  By the convexity of the $L_1$-norm, $\lonenorm{\hat{\bbeta}^t}-\lonenorm{\hat{\bbeta}}
  - [\nabla \lonenorm{\hat{\bbeta}}]^T(\hat{\bbeta}^t-\hat{\bbeta}) \geq 0$. Hence,
  \begin{equation}
    \label{eq7.24}
    \phi(\hat{\bbeta}^t)-\phi(\hat{\bbeta})\geq\frac{\gamma_l}{2}
    \ltwonorm{\hat{\bbeta}^t-\hat{\bbeta}}^2
    -\tau_l \lonenorm{\hat{\bbeta}^t-\hat{\bbeta}}^2.
  \end{equation}
  Next, we bound $\lonenorm{\hat{\bbeta}^t-\hat{\bbeta}}$. It follows from Lemma 3 of
  \cite{ANW12} that
  \begin{equation*}
    \lonenorm{\hat{\bbeta}^t-\hat{\bbeta}}\leq 2 \left ( 2 \sqrt{S_{\alpha
          \eta}}\ltwonorm{\hat{\bbeta}^t-\hat{\bbeta}} + 4 \sqrt{|S_{\alpha
          \eta}|}\ltwonorm{\hat{\bbeta}-\betaa}+4
      \lonenorm{\betas_{S_{\alpha\eta}^c}}+ {\delta^2}/{\lambda_n} \right ),
  \end{equation*}
  where $\delta$ is defined as in (a). Then, by the Cauchy-Schwartz inequality,
  \begin{equation}
    \label{eq7.25}
    \lonenorm{\hat{\bbeta}^t-\hat{\bbeta}}^2\leq 16 \left (4
      |S_{\alpha\eta}|\ltwonorm{\hat{\bbeta}^t-\hat{\bbeta}}^2 + 16 |S_{\alpha \eta}|
      \ltwonorm{\hat{\bbeta}-\betas_{\alpha}}^2 + 16 \lonenorm{\betas_{S_{\alpha \eta}^c}}^2 +
      {\delta^4}/{\lambda_n^2} \right ).
  \end{equation}
  Equations (\ref{eq7.24}) and (\ref{eq7.25}) together with results in (a) imply that,
  \begin{equation*}
    \delta^2\geq \frac{\gamma_l}{2} \ltwonorm{\hat{\bbeta}^t-\hat{\bbeta}}^2-16 \tau_l\left(4
      |S_{\alpha\eta}|\ltwonorm{\hat{\bbeta}^t-\hat{\bbeta}}^2 + 16 |S_{\alpha \eta}|
      \ltwonorm{\hat{\bbeta}-\betas_{\alpha}}^2 + 16 \lonenorm{\betas_{S_{\alpha \eta}^c}}^2 +
      {\delta^4}/{\lambda_n^2}\right).
  \end{equation*}
  Letting $\tilde{\gamma}_l= {\gamma_l}/{2}-64\tau_l|S_{\alpha \eta}|$, we have
  \begin{equation}
    \label{eq7.26}
    \ltwonorm{\hat{\bbeta}^t-\hat{\bbeta}}^2\leq
    \frac{1}{\tilde{\gamma}_l}\left(\delta^2+\frac{16\tau_l\delta^4}{\lambda_n^2}\right)
    +\frac{256\tau_l}{\tilde{\gamma}_l}(|S_{\alpha\eta}|\ltwonorm{\hat{\bbeta}-\betaa}^2
    +\lonenorm{\betas_{S_{\alpha\eta}^c}}^2).
  \end{equation}
  We now bound the second term in (\ref{eq7.26}).  By (\ref{eq7.16}) and (\ref{eq7.17}), we have
  \begin{equation}
    \label{eq7.27}
    \begin{split}
      |S_{\alpha\eta}|\ltwonorm{\hat{\bbeta}-\betaa}^2+\lonenorm{\betas_{S_{\alpha\eta}^c}}^2&\leq
      \Rq\eta^{-q} \ltwonorm{\hat{\bbeta}-\betaa}^2 +\Rq^2\eta^{2-2q}\\&\leq
      \Rq\kappa_{\lambda}^{-q}\left(\frac{\log p}{n}\right)^{-q/2}\ltwonorm{\hat{\bbeta}-\betaa}^2
      +\kappa_{\lambda}^{-q}\Rq^2\left(\frac{\log p}{n}\right)^{1-q}\\
      &\leq \kappa_{\lambda}^{-q}\Rq\left(\frac{\log
          p}{n}\right)^{-q/2}\left[\ltwonorm{\hat{\bbeta}-\betaa}^2 +\Rq\left(\frac{\log
            p}{n}\right)^{1-(q/2)}\right].
    \end{split}
  \end{equation}
  Meanwhile, from (a) we have
  \begin{equation}
    \label{eq7.28}
    \begin{split}
      \delta^2&=\frac{\varepsilon^2}{1-\kappa}=\frac{8\xi\chi}{1-\kappa}\left(6 \sqrt{|S_{\alpha\eta}|}
        \ltwonorm{\hat{\bbeta}-\betaa}+8 \lonenorm{\betas_{S_{\alpha\eta}^c}}\right)^2 \\
      &\leq \frac{8\xi\chi}{1-\kappa}(72|S_{\alpha\eta}|\ltwonorm{\hat{\bbeta}-\betaa}^2+128
      \lonenorm{\betas_{S_{\alpha\eta}^c}}^2)\\
      &\leq \frac{1024 \xi\chi}{1-\kappa}(|S_{\alpha\eta}|
      \ltwonorm{\hat{\bbeta}-\betaa}^2+\lonenorm{\betas_{S_{\alpha\eta}^c}}^2).
    \end{split}
  \end{equation}
  Since $\bar{\gamma}_l\asymp 1$, $\kappa\asymp 1$, $\xi\asymp 1$, $\chi\asymp \frac{\log p}{n}$,
  and $\tau_l\asymp \frac{\log p}{n}$, it follows from (\ref{eq7.26}), (\ref{eq7.27}) and
  (\ref{eq7.28}) that
  \begin{equation*}
    \ltwonorm{\hat{\bbeta}^t-\hat{\bbeta}}^2=O \left(R_q\left(\frac{\log p}{n}\right)^{1-(q/2)}
      \left[\ltwonorm{\hat{\bbeta}-\betaa}^2+R_q\left(\frac{\log p}{n}\right)^{1-(q/2)}\right] \right) .
  \end{equation*}
\end{proof}

\begin{proof}[\textbf{Proof of Theorem \ref{thm5}.}]
  The proof follows the same spirit of the proof of Proposition 2.4 of \cite{CAT12}. The influence
  function $\psi(x)$ satisfies
  \begin{equation*}
    -\log(1-x+x^2)\leq \psi(x)\leq \log(1+x+x^2).
  \end{equation*}
  Using this and independence, with $r(\theta)=\frac{1}{\alpha n}\sum_{i=1}^n \psi[\alpha(Y_i-\theta)]$, we have
  \begin{align*}
    \E \left\{\exp[\alpha n r(\theta)] \right\}
    & \leq \Bigl ( \E \{\exp \{ \psi[\alpha(Y_i-\theta)] \} \Big )^n \\
    & \leq \{1+\alpha(\mu-\theta)+\alpha^2[\sigma^2+(\mu-\theta)^2] \}^n\\
    &\leq \exp \left\{n\alpha(\mu-\theta)+n\alpha^2[ v^2+(\mu-\theta)^2] \right\}.
  \end{align*}
  Similarly, $\E \left\{\exp[-\alpha n r(\theta)] \right\}\leq \exp
  \left\{-n\alpha(\mu-\theta)+n\alpha^2[ v^2+(\mu-\theta)^2] \right\}$. Define
  \begin{align*}
    B_+(\theta)&=\mu-\theta+\alpha[ v^2+(\mu-\theta)^2]+\frac{\log (1/\delta)}{n\alpha}\\
    B_-(\theta)&=\mu-\theta-\alpha[ v^2+(\mu-\theta)^2]-\frac{\log (1/\delta)}{n\alpha}
  \end{align*}
  By Chebyshev inequality,
  \begin{equation*}
    P(r(\theta)>B_+(\theta))\leq \frac{\E \left\{\exp[\alpha n r(\theta)] \right\}}
    {\exp\{\alpha n(\mu-\theta)+n\alpha^2[ v^2+(\mu-\theta)^2]+\log(1/\delta) \}}\leq \delta
  \end{equation*}
  Similarly, $P(r(\theta)<B_-(\theta))\leq \delta$.

  Let $\theta_+$ be the smallest solution of the quadratic equation $B_+(\theta_+)=0$ and $\theta_-$
  be the largest solution of the equation $B_-(\theta_-)=0$. Under the assumption that $\frac{\log
    (1/\delta)}{n}\leq 1/8$ and the choice of $\alpha=\sqrt{\frac{\log
      (1/\delta)}{n v^2}}$, we have $\alpha^2 v^2+\frac{\log (1/\delta)}{n}\leq
  1/4$. Therefore,
  \begin{align*}
    \theta_+&=\mu+2\left(\alpha v^2+\frac{\log(1/\delta)}{\alpha n}\right) \left(1+ \sqrt{1-4
        \left(\alpha^2 v^2+\frac{\log (1/\delta)}{n} \right)} \right)^{-1}\\
    &\leq \mu+2\left(\alpha v^2+\frac{\log(1/\delta)}{\alpha n}\right).
  \end{align*}
  Similarly,
  \begin{align*}
    \theta_-&=\mu-2\left(\alpha v^2+\frac{\log(1/\delta)}{\alpha n}\right) \left(1+ \sqrt{1-4
        \left(\alpha^2 v^2+\frac{\log (1/\delta)}{n} \right)} \right)^{-1}\\
    &\geq \mu-2\left(\alpha v^2+\frac{\log(1/\delta)}{\alpha n}\right).
  \end{align*}
  With $\alpha=\sqrt{\frac{\log(1/\delta)}{n v^2}}$, $\theta_+\leq \mu+4 v
  \sqrt{\frac{\log (1/\delta)}{n}}$, $\theta_-\geq \mu-4 v \sqrt{\frac{\log
      (1/\delta)}{n}}$. Since the map $\theta \mapsto r(\theta)$ is non-increasing, under event
  $\{B_-(\theta)\leq r(\theta)\leq B_+(\theta) \}$
  \begin{equation*}
    \mu-4 v \sqrt{\frac{\log(1/\delta)}{n}} \leq \theta_-\leq \hat{\mu}_{\alpha}
    \leq \theta_+\leq \mu+4 v \sqrt{\frac{\log(1/\delta)}{n}},
  \end{equation*}
  i.e. $|\hat{\mu}_{\alpha}-\mu|\leq 4 v \sqrt{\frac{\log (1/\delta)}{n}}$. Meanwhile,
  $P(B_-(\theta)\leq r(\theta)\leq B_+(\theta))>1-2\delta$.
\end{proof}

\begin{proof}[\textbf{Proof of Theorem \ref{thm6}.}]
  First, we prove that the approximation error rate $\ltwonorm{\betacs-\betas}=O(\alpha^{k-1})$,
  where $\betacs=\argmin_{\bbeta} \E \ell_{\alpha}^c(y-\bx^T\betas)$ is the population
  minimizer under the Catoni loss. Let $g_{\alpha}(x)=\ell(x)-\ell_{\alpha}^c(x)=\int_0^x
  [2t-\frac{2}{\alpha}\psi_c(\alpha t)] \text{d}t$. It follows from (\ref{eq7.2}) that
  \begin{equation*}
    \E[\ell(y-\bx^T\betaa)-\ell(y-\bx^T\betas)]\leq
    \E[|g_{\alpha}'(y-\bx^T\tilde{\bbeta})\bx^T(\betacs-\betas)|],
  \end{equation*}
  where $\tilde{\bbeta}$ is a vector lying between $\betas$ and $\betacs$. Since $|(\psi_c)'''|\leq
  3$, by the second-order Taylor expansion with an integral remainder,
  \begin{equation}
    \label{eq7.29}
    |g_{\alpha}'(x)|=|2x-\frac{2}{\alpha}\psi_c(\alpha x)|=\left|\frac{\alpha^2}{3}\int_0^x (\psi_c)'''
      (\alpha s)(x-s)^2 ds\right|
    \leq \alpha^2|x|^3.
  \end{equation}
  Hence, by the Cauchy-Schwartz inequality,
  \begin{align*}
    \E[\ell(y-\bx^T\bbeta_{\alpha}^{c\ast})-\ell(y-\bx^T\betas)]&\leq \alpha^2
    \E[|y-\bx^T\tilde{\bbeta}|^3|\bx^T(\bbeta_{\alpha}^{c\ast}-\betas)|] \\
    &\leq 4\alpha^2\E[(M_3+|\bx^T(\tilde{\bbeta}-\betas)|^3)|\bx^T(\bbeta_{\alpha}^{c\ast}-\betas)|]\\
    &\leq 4\alpha^2[\lambda_{\max}(\E[(M_3+|\bx^T(\tilde{\bbeta}-\betas)|^3)^2
    \bx\bx^T])]^{1/2} \ltwonorm{\bbeta_{\alpha}^{c\ast}-\betas},
  \end{align*}
  which is of order $O(\alpha^2 \ltwonorm{\bbeta_{\alpha}^{c\ast}-\betas})$, as
  $\lambda_{\max}(\E[(M_3+|\bx^T(\tilde{\bbeta}-\betas)|^3)^2 \bx\bx^T])=O(1)$ by an analogous
  argument as in the proof of Theorem \ref{thm1}. Similarly as in (\ref{eq7.1}),
  \begin{equation*}
    \E[\ell(y-\bx^T\bbeta_{\alpha}^{c\ast})-\ell(y-\bx^T\betas)]\geq \kappa_l
    \ltwonorm{\betaa-\betas}^2.
  \end{equation*}
  Hence, $\ltwonorm{\betacs-\betas}=O(\alpha^2)$. If $\epsilon$ only has the second moment exist,
  by a first-order Taylor expansion of $g_{\alpha}'(x)$ similarly as in (\ref{eq7.29}), we have
  $\ltwonorm{\bbeta_{\alpha}^{c\ast}-\betas}=O(\alpha)$. Next, since $(\psi_c)'(0)=1$, by the same
  argument as in the proof of Lemma \ref{lem3}, RSC holds for Catoni's loss with probability no
  less than $1-c_1\exp(-c_2n)$. Hence, similarly as in Theorem \ref{thm2},
  $\ltwonorm{\hat{\bbeta}-\betacs}=O(\sqrt{\Rq}[(\log p)/n]^{1/2-q/4})$. This together with
  $\ltwonorm{\betacs-\betas}=O(\alpha^{k-1})$ completes the proof.
\end{proof}

\begin{proof}[\textbf{Proof of Theorem \ref{thm7}.}]
  First of all, observe that
  \begin{align*}
    \hat{\sigma}^2-\sigma^2 &=\frac{1}{K} \sum_{k=1}^K \frac{1}{m} \sum_{i \in \text{fold } k}
    \left(\epsilon_i-(\bx_i^T\hat{\bbeta}^{(-k)}-\bx_i^T\betas)\right)^2-\sigma^2\\
    &=\frac{1}{n}\sum_{i=1}^n\epsilon_i^2-\sigma^2-\frac{1}{K}\sum_{k=1}^K \frac{2}{m}\sum_{i \in
      \text{fold } k} \epsilon_i\bx_i^T\left(\hat{\bbeta}^{(-k)}-\betas\right)
    +\frac{1}{K}\sum_{k=1}^K \frac{1}{m} \sum_{i \in \text{fold } k}
    \left(\bx_i^T(\hat{\bbeta}^{(-k)}-\betas)\right)^2.
  \end{align*}
  Given that $\E\epsilon^4$ exists, by Central Limit Theorem, $\sqrt{n}(\frac{1}{n}
  \sum_{i=1}^n\epsilon_i^2-\sigma^2)\stackrel{\mathcal{D}}{\to} \mathcal{N}(0,
  \E\epsilon^4-\sigma^4)$. Let $z_i=\bx_i^T(\hat{\bbeta}^{(-k)}-\betas)$.
  We now need to prove that the last two terms are negligible.   Conditioning on data outside the $k$th fold, 
  \begin{equation*}
   \E\left\{\frac{1}{m} \left(\sum_{i\in \text{fold }k} \epsilon_i z_i \right)^2
    \right\} = \E(\epsilon_i^2z_i^2)\leq \sigma^2\kappa_0^2
    \ltwonorm{\hat{\bbeta}^{(-k)}-\betas}^2.
  \end{equation*}
  Hence, $m^{-1/2} \sum_{i\in \text{fold }k} \epsilon_i
  \bx_i^T(\hat{\bbeta}^{(-k)}-\betas)=O_P\left(\ltwonorm{\hat{\bbeta}^{(-k)}-\betas}\right)=o_P(1)$,
  where the last equality follows from Theorem \ref{thm3}. By an analogous argument, we have
  \begin{align*}
    \frac{1}{m} \sum_{i \in \text{fold } k} \left(\bx_i^T(\hat{\bbeta}^{(-k)}-\betas)\right)^2 {=}
    O_p\left(\ltwonorm{\hat{\bbeta}^{(-k)}-\betas}^2\right)=O_p(\max\{\alpha^{2(k-1)},{\Rq}[(\log
    p)/n]^{1-q/2} \})=o(1/\sqrt{n}).
  \end{align*}
  This completes the proof of the Theorem.
\end{proof}

\end{document}